\documentclass{article}

\usepackage{amssymb}
\usepackage{amsfonts}
\usepackage{amsmath}
\usepackage{url}
\usepackage{latexsym}

\newtheorem{theorem}{Theorem}[section]
\newtheorem{corollary}[theorem]{Corollary}
\newtheorem{lemma}[theorem]{Lemma}
\newtheorem{definition}[theorem]{Definition}

\newtheorem{rem}[theorem]{Remark}
\newtheorem{example}[theorem]{Example}

\newtheorem{problem}{Problem}

\def\deg{\mathop{\rm deg }\nolimits}
\def\rank{\mathop{\rm rank}\nolimits}
\def\diag{\mathop{\rm diag }\nolimits}
\def\lcm{\mathop{\rm lcm }\nolimits}

\newcommand{\efe}{\mathbb F}
\newcommand{\FF}{\mathbb F}
\newcommand{\CC}{\mathbb C}
\newcommand{\RR}{\mathbb R}

\newcommand{\la}{s}

\newcommand{\ba}{{\bf a}}
\newcommand{\bb}{{\bf b}}
\newcommand{\bc}{{\bf c}}
\newcommand{\bd}{{\bf d}}
\newcommand{\bg}{{\bf g}}
\newcommand{\bu}{{\bf u}}
\newcommand{\bv}{{\bf v}}

\title{Row completion of polynomial and rational matrices II
  \thanks{
This work was 
supported by grant PID2021-124827NB-I00 funded by MCIN/AEI/ 10.13039/501100011033 and by ``ERDF A way of making Europe'' by the ``European Union''.
The first and third authors were also supported 
by grant GIU21/020 funded by UPV/EHU.
  }
}

\author{Agurtzane Amparan\thanks{Departamento de Matem\'aticas, Universidad del Pa\'is Vasco UPV/EHU, Bilbao, Spain, {agurtzane.amparan@ehu.eus}, {silvia.marcaida@ehu.eus}}
\and Itziar Baraga\~na\thanks{Departamento de Ciencia de la Computaci\'on e I.A., 
Universidad del Pa\'{\i}s Vasco UPV/EHU, Donostia-San Sebasti\'an, 
Spain, {itziar.baragana@ehu.eus}}
\and Silvia Marcaida\footnotemark[2]
\and Alicia Roca\thanks{Departamento de Matem\'atica Aplicada, IMM, Universitat Polit\`ecnica de Val\`encia, Valencia, Spain,   {aroca@mat.upv.es}}}

\date{}

\begin{document}

\maketitle

\begin{abstract}
We study the row completion problem of polynomial and rational matrices with partial prescription of the structural data. 
The prescription of the complete structural data 
 has been solved 
 in \cite{AmBaMaRo25}, where several results of prescription of some of the four types of invariants  composing  the structural data
have also been obtained. In this paper we  deal with the cases not analyzed there.
More precisely, we solve the row completion problem of a
rational or a polynomial  matrix when we prescribe the infinite (finite)
structure and the column and/or the row minimal indices, and when only
the column and/or row minimal indices are prescribed. The necessity of the
conditions obtained are valid over arbitrary
ﬁelds, but in some cases the proof of the suﬃciency requires working 
over algebraically closed fields. By transposition the results obtained hold for the corresponding column completion problems.
\end{abstract}


{\bf AMS:}
   15A18, 15A54, 15A83

\section{Introduction}

The problem addressed in this paper consists in characterizing the existence of a rational or a polynomial matrix  with some rows (or columns) and part of the structural data  prescribed. 
Formally stated, the problem is the following:

\begin{problem}\label{problemrat}
Let $R(s)\in\efe(s)^{m\times n}$ be a rational (polynomial) matrix. Find necessary and sufficient conditions for the existence of a rational (polynomial) matrix $\widetilde{W}(s)\in \efe(s)^{z\times n}$
such that $\begin{bmatrix}R(s)\\\widetilde{W}(s)\end{bmatrix}$  has  part of the structural data prescribed.
\end{problem}

The {\em complete structural data} of a rational matrix are formed by the invariant rational functions (also known as the finite structure), the invariant orders at $\infty$ (known as the  infinite structure), and the column and row minimal indices (the singular structure) (\cite{AnDoHoMa19}). We remark that for a polynomial matrix the  invariant rational functions are the invariant factors, and the smallest  invariant order at $\infty$ is minus the degree of the  matrix.
Therefore, when we prescribe the invariant orders at $\infty$ of a polynomial matrix, we also prescribe its degree. 

If a rational matrix is transposed, the finite and infinite structures remain the same, and the column and row minimal indices are interchanged. As a consequence, by transposition, the row completion results obtained in this paper  also hold for the corresponding column completion.

The prescription of  only part of the  structural data allows us to characterize when  a rational  matrix can be completed to achieve certain  invariants, leaving some freedom for the rest of them. 
Moreover, in some cases, it permits us to characterize  a  completion under  partial knowledge of the structural data of the  submatrix. Since the complete structural data   consist of four types of invariants, the analysis of the 
whole and
partial prescription of all possible cases leads to  15 different problems.

This paper is a  continuation of  \cite{AmBaMaRo24},  \cite{AmBaMaRo24_2}, and \cite{AmBaMaRo25}.
In  \cite{AmBaMaRo24} and \cite{AmBaMaRo24_2} we characterize the complete structural data or part of it of a polynomial matrix when some of its rows are prescribed
 (this is the {\em polynomial case of Problem \ref{problemrat}}) and, additionally, the degree of the completed matrix coincides with that of the submatrix. 
  In \cite{AmBaMaRo24} we provide a solution to  6 of the 15 resulting   problems; specifically, the  prescription of the  complete structural data, of  the complete structure but the row (column) minimal indices, and the finite and/or infinite  structure. In \cite{AmBaMaRo24_2}  the remaining cases are solved. 

The next step was to remove the restriction on the degree of the completed matrix, allowing the degree to increase, therefore  achieving more possibilities of completion. 
For instance, in Example \ref{exalgclosed} below it is possible to  obtain the desired invariants  if the degree of the completed matrix is allowed to be greater than that of the prescribed submatrix, but it is not possible otherwise. This idea led to the polynomial case of Problem \ref{problemrat}. Finally, we extended the problem to rational matrices as stated in Problem \ref{problemrat}.
 In \cite{AmBaMaRo25} we solved it in the cases analogous to those solved in \cite{AmBaMaRo24}. In this paper we complete the research solving the remaining cases.

 As we show along the paper,  when the finite structure is prescribed, a solution to the polynomial case can be derived from that of the  rational one, but in general the characterizations  obtained for the polynomial case  are not particular cases of the corresponding rational completion results.
 
The matrix completion problem has proven to be an important problem,  from the mathematical  and applied  points of view. 
 It has been addressed both theoretically and numerically, generating a really extensive  literature. Covering it is an unattainable task in this paper. We would like to mention  just  a few results, relevant for our approach to the problem. 
It has been studied for square constant matrices prescribing similarity invariants \cite{CaSi92,  Ol69, FaLe59,   Sa79, Mi58,   Si87, Th79, ZaLAA87}, for rectangular matrices and feedback equivalence invariants \cite{BaZa90, Do05, ZaLAA88}, for matrix pencils and Kronecker invariants  \cite{Ba89, CaSi91, Do08, Do10, DoSt19, FuSi99, LoMoZaZaLAA98}, for polynomial matrices and unimodular equivalence   \cite{Sa79, Th79}, and finally for polynomial and rational matrices and structural data \cite{AmBaMaRo24,AmBaMaRo24_2, AmBaMaRo25}. See also the references therein.

The paper is organized as follows. Section \ref{secprelimin} contains the notation, definitions and previous results. In particular, we restate the result in 
\cite[Theorem 3.11]{AmBaMaRo25}, where we present  a solution to the row completion  problem  for rational matrices when the complete structural data are prescribed (see Theorem \ref{prescr4rat}  below). 
Section \ref{secremaining} is devoted to solve the row completion problem for rational and polynomial matrices in the cases not solved in \cite{AmBaMaRo25}. 
In Subsection \ref{subsecinfsing} we deal with the problems of prescribing the infinite and singular structures. Different cases are solved: prescription of the infinite and singular structures for polynomial (rational) matrices in Theorem \ref{theoprescrioirmicmipol} (Theorem \ref{theoprescrioirmicmirat}), prescription of the  infinite structure and the column minimal indices  for polynomial (rational) matrices in Theorem \ref{corprescrioicmipol} (Theorem \ref{corprescrioicmirat}), and
 prescription of the infinite structure and the row minimal indices  for polynomial (rational) matrices in Theorem \ref{corprescrioirmipol} (Theorem \ref{corprescrioirmirat}).
In Subsection \ref{subsecfinsing} we address the prescription of the finite and singular structures. As mentioned, when the finite structure is prescribed, the polynomial case  follows from the rational case. The prescription of the finite and singular structures, of the finite structure and the column minimal indices, and of the  finite structure and the row minimal indices are solved in Theorems \ref{theoprescrifrmicmirat}, \ref{corprescrifcmi} and  \ref{corprescrifrmi}, respectively. Subsection \ref{subsecsing} deals with  the prescription of the singular structure. The conditions obtained for the polynomial case are the same as those obtained for the rational case. The prescription of the singular structure, of the row minimal indices, and of the column minimal indices are analyzed in 
Theorems \ref{theoprescrrmicmirat}, \ref{corprescrrmi} and \ref{corprescrcmi}, respectively.
 We must mention that the proofs of Theorems \ref{theoprescrioirmicmipol}, \ref{corprescrioicmipol} and \ref{corprescrioirmipol} are very similar to those of the corresponding theorems  in \cite{AmBaMaRo24_2}.  As they are rather intricate, we have included them in \ref{secappendix} for the reader's convenience.

\section{Preliminaries}\label{secprelimin}

Let $\FF$ be a field. 
The ring of polynomials in the indeterminate $s$ with coefficients in $\FF$
is denoted by $\FF[s]$,  $\FF(s)$ is the field of fractions of
$\FF[s]$, i.e., the field of rational functions over $\FF$, and $\FF_{pr}(s)$ is the ring of proper rational functions, i.e., the rational functions with degree of the denominator at least the degree of the numerator. 
A polynomial in $\FF[s]$ is \textit{monic} if its leading coefficient is 1. 
  Given two polynomials $\alpha(s), \beta(s)$, by $\alpha(s)\mid \beta(s)$ we mean that $\alpha(s)$ is a divisor of $\beta(s)$, by $\lcm(\alpha,\beta)$, the monic least common multiple of $\alpha(s)$ and $\beta(s)$, and by $\gcd(\alpha,\beta)$, the monic greatest common divisor of $\alpha(s)$ and $\beta(s)$.

   In this work we deal with  
    polynomial chains $\alpha_1(s)\mid \dots \mid \alpha_r(s)$, where $\alpha_i(s)\in \FF[s]$, and take $\alpha_1(s)=1$ for $i<1$ and $\alpha_i(s)=0$ for $i>r$.
If the polynomial chain is ordered as $\varphi_r(s)\mid \dots \mid \varphi_1(s)$,  we take $\varphi_i(s)=1$ for $i>r$ and $\varphi_i(s)=0$ for $i<1$.

Along the paper, if $a_1\geq\dots\geq a_r$ is a  decreasing sequence of integers, we write $\ba=(a_1, \dots, a_r)$. If $a_r\geq 0$,
the sequence is called a 
\textit{partition}. When necessary, we take $a_i=+\infty$ for $i<1$ and 
$a_i=-\infty$ for $i>r$.  If $b_1\leq\dots\leq b_r$ is an increasing sequence of integers,  we take $b_i=-\infty$ for $i<1$ and 
$b_i=+\infty$ for $i>r$.

We denote by  $\FF^{m\times n}$, $\FF[s]^{m\times n}$, $\FF(s)^{m \times n}$, and $\FF_{pr}(s)^{m \times n}$ the sets   of $m \times n$ matrices with elements in $\FF$, $\FF[s]$, $\FF(s)$, and $\FF_{pr}(s)$, respectively.
The \textit{degree} of a polynomial matrix $P(s)$, $\deg(P(s))$, is the highest degree of the entries of $P(s)$.
A matrix $U(s)\in\efe[s]^{n\times n}$ is said \textit{unimodular} if it has inverse in $\efe[s]^{n\times n}$, while a matrix $B(s)\in\efe_{pr}(s)^{n\times n}$ is said \textit{biproper} if it has inverse in $\efe_{pr}(s)^{n\times n}$.

Two rational matrices $R_1(s), R_2(s)\in \FF(s)^{m\times n}$ are \textit{unimodularly equivalent} if there exist unimodular matrices $U_1(s)\in\efe[s]^{m\times m}$ and $U_2(s)\in\efe[s]^{n\times n}$ such that $R_2(s)=U_1(s)R_1(s)U_2(s)$.  Let $R(s)\in \FF(s)^{m\times n}$ be a rational matrix  of $\rank(R(s))=r$. A canonical form for the unimodular equivalence of $R(s)$   is the \textit{Smith--McMillan form}
$$
\begin{bmatrix}\diag\left(\frac{\eta_1(s)}{\varphi_1(s)},\dots,\frac{\eta_r(s)}{\varphi_r(s)}\right)&0\\0&0\end{bmatrix},$$
where $\eta_1(s)\mid  \dots \mid \eta_r(s)$ and $\varphi_r(s)\mid  \dots \mid \varphi_1(s)$ are monic polynomials, and $\frac{\eta_1(s)}{\varphi_1(s)},\dots,\frac{\eta_r(s)}{\varphi_r(s)}$ are irreducible rational functions
 known as the \textit{invariant rational functions} of $R(s)$. We also refer to them as the \textit{finite structure} of $R(s)$. The polynomial $\varphi_1(s)$ is the monic least common denominator of the entries of $R(s)$ (see, for instance, \cite[Chapter 3, Section 4]{Rose70}). 

If the rational matrix is a polynomial matrix $P(s)$,  then $\varphi_1(s)=\dots=\varphi_r(s)=1$, the polynomials $\eta_1(s)\mid\dots\mid\eta_r(s)$ are the \textit{invariant factors} of $P(s)$, and the Smith-McMillan form is its \textit{Smith normal form} (\cite[Chapter 1, Section 1]{Rose70}).
Recall that (see, for instance, \cite[p. 261]{LaTi85})
$$
\eta_1(s)\cdots\eta_k(s)=\gcd\{m_k(s):  \mbox{$m_k(s) = $ minor of order  $k$ of $P(s)$}\},\quad 1\leq k\leq r.
$$

Two rational matrices $R_1(s), R_2(s)\in \FF(s)^{m\times n}$ are \textit{equivalent at infinity} if there exist biproper matrices $B_1(s)\in\efe_{pr}(s)^{m\times m}$ and $B_2(s)\in\efe_{pr}(s)^{n\times n}$ such that $R_2(s)=B_1(s)R_1(s)B_2(s)$. Let $R(s)\in \FF(s)^{m\times n}$ be a rational matrix  of rank $r$. A canonical form for the equivalence at infinity of $R(s)$ is the  \textit{Smith--McMillan  form at infinity}  
$$ 
\begin{bmatrix}\diag\left(s^{-\tilde p_1},\ldots,s^{-\tilde p_r}\right)&0\\0&0\end{bmatrix},
$$
where  $\tilde p_1\leq\cdots\leq\tilde p_r$ are integers called the \textit{invariant orders at infinity} of $R(s)$  (see, for instance, \cite{Vard91}).
In \cite{AnDoHoMa19}, the sequence  of invariant orders at $\infty$  is called  the \textit{structural index sequence of $R(s)$ at $\infty$.} For a polynomial matrix $P(s)$, it  was proved in \cite[p. 102]{Vard91} that 
$$\sum_{i=1}^k\tilde{p}_i=-\max\{\deg(m_k(s)):  \mbox{$m_k(s) = $ minor of order  $k$ of  $P(s)$}\},\quad 1\leq k\leq r.
$$
In particular, $\tilde{p}_1=-\deg(P(s))$.

We recall now the \textit{singular structure} of a rational matrix. Denote by $\mathcal{N}_\ell (R(\la))$ and $\mathcal{N}_r (R(\la))$ the \textit{left} and
\textit{right null-spaces} over $\FF(\la)$ of $R(\la)$, respectively, i.e.,
if $R(\la)\in\FF(\la)^{m\times n}$,
\[
\begin{array}{l}
\mathcal{N}_\ell (R(\la))=\{x(\la)\in\FF(\la)^{m \times 1}: x(\la)^TR(\la)=0\},\\
\mathcal{N}_r (R(\la))=\{x(\la)\in\FF(\la)^{n \times 1}: R(\la)x(\la)=0\},
\end{array}
\]
which are vector subspaces
of $\FF(\la)^{m \times 1}$ and $\FF(\la)^{n \times 1}$, respectively. Given  a subspace $\mathcal{V}$ of
$\FF(\la)^{m \times 1}$ it is possible to find a basis consisting of vector
polynomials; it is enough to  take an arbitrary basis and multiply each vector
by a  least common multiple of the denominators of its  entries.
The \textit{order} of a polynomial basis is defined as the sum of
the degrees of its vectors  (see \cite{Fo75}). A \textit{minimal basis} of  $\mathcal{V}$
is  a polynomial basis  with least order among the
polynomial bases of $\mathcal{V}$. 
The  degrees of the vector polynomials  of a minimal basis, increasingly ordered,  are  always the same (see \cite{Fo75}), and are called the \textit{minimal indices} of $\mathcal{V}$.

 A   \textit{right  (left) minimal basis} of a rational matrix $R(\la)$ is a minimal basis of $\mathcal{N}_r (R(\la))$ ($\mathcal{N}_\ell (R(\la))$).
The \textit{right  (left) minimal indices} of $R(\la)$ are the minimal indices of
$\mathcal{N}_r (R(\la))$ ($\mathcal{N}_\ell (R(\la))$). 
From now on in this paper, we  work with the right (left) minimal indices decreasingly ordered, and we  refer to them as the \textit{column (row) minimal indices} of $R(s)$.  Notice that a rational matrix $R(s)\in \FF(s)^{m\times n}$ of $\rank(R(s))=r$ has $m-r$ row and $n-r$ column minimal indices.

Given a rational matrix $R(s)\in \FF(s)^{m\times n}$
of $\rank(R(s))=r$, the \textit{complete structural data} consist of four components
(see \cite[Definition 2.15]{AnDoHoMa19}): the invariant rational functions $\frac{\eta_1(s)}{\varphi_1(s)},\dots,   \frac{\eta_r(s)}{\varphi_r(s)}$,
the invariant orders at infinity  $\tilde p_1\leq \cdots \leq  \tilde p_{r}$, the row minimal indices $(u_1, \dots, u_{m-r})$ and the column minimal indices $(c_1, \dots, c_{n-r})$.
Observe that the complete structural data of a rational matrix determine its rank.

Given two integers $n$ and $m$, whenever $n>m$ we take  $\sum_{i=n}^{m}=0$. In the same way, if a condition is stated for $n\leq i\leq m$ with $n>m$, we understand that the condition 
 is trivially satisfied.


Let   $\ba= (a_1,  \ldots, a_m)$ and $\bb= (b_1, \ldots, b_m)$ be  two sequences of integers. It is said that   $\ba$ is \textit{majorized} by $\bb$ (denoted by $\ba \prec \bb$) if $\sum_{i=1}^k a_i \leq \sum_{i=1}^k b_i $ for $1 \leq k \leq m-1$ and $\sum_{i=1}^m a_i =\sum_{i=1}^m b_i$ (this is an extension to sequences of integers of the definition of majorization given for partitions in \cite{HLP88}).
We introduce next the concept of \textit{generalized majorization}. 
\begin{definition}{\rm \cite[Definition 2]{DoStEJC10}}
Let $\bd = (d_1, \dots, d_m)$, $\ba=(a_1, \dots, a_s)$ and $\bg=(g_1, \dots, g_{m+s})$  be sequences of  integers.
We say that  $\bg$ is majorized by $\bd$ and $\ba$  $(\bg \prec' (\bd,\ba))$ if
\begin{equation}\label{gmaj1}
d_i\geq g_{i+s}, \quad 1\leq i\leq m,
\end{equation}
\begin{equation}\label{gmaj2}
\sum_{i=1}^{h_j}g_i-\sum_{i=1}^{h_j-j}d_i\leq \sum_{i=1}^j a_i, \quad 1\leq j\leq s,
\end{equation}
where $h_j=\min\{i:  d_{i-j+1}<g_i\}$, $1\leq j\leq s$  $(d_{m+1}=-\infty)$,
\begin{equation}\label{gmaj3}
\sum_{i=1}^{m+s}g_i=\sum_{i=1}^md_i+\sum_{i=1}^sa_i.
\end{equation}
\end{definition}
\begin{rem}\label{gmajcp} 
Recall that $d_i=+\infty$ for $i\leq 0$ and $d_i=-\infty$ for $i>m$. Therefore $j\leq h_j\leq m+j$. If $s=0$,  then $\bg \prec' (\bd,\ba)$ is equivalent to $\bd=\bg$,
 and if $m=0$ it reduces to  $\bg \prec \ba$.   Also, if $\bg \prec' (\bd,\ba)$ and $\ba\prec \hat \ba$ for some $\hat \ba$, then  $\bg \prec' (\bd,\hat \ba)$.
\end{rem}
 
Given two sequences of integers $\bu = (u_1, \dots, u_{p})$ and $\bb = (b_1, \dots, b_{y})$ the union, $\bu\cup \bb$,  is the decreasingly ordered sequence  of the $p+y$ integers of $\bu$ and $\bb$.
The next lemma is satisfied.
\begin{lemma} {\rm (\cite[Lemma 4.4]{AmBaMaRo24})}\label{lemmacup}
 Let $\bu = (u_1, \dots, u_{p})$ and $\bb = (b_1, \dots, b_{y})$    be sequences of  integers.  
Then
$
\bu\cup \bb\prec'(\bu, \bb).
$
\end{lemma}

The following technical lemma will be used later.

\begin{lemma}\label{lemmahx}
Let $x, r, n$ be integers such that $0\leq x \leq n-r$, and 
let $\bd=(d_1, \dots, d_{n-r-x})$, $\bc=(c_1, \dots, c_{n-r})$, $\ba=(a_1, \dots, a_x)$ be partitions. 
Let $h_x=\min\{i: d_{i-x+1}<c_i\}$.
\begin{enumerate}
\item
If $\bc\prec'(\bd, \ba)$, then
\begin{equation}\label{eqdc}
d_i=c_{i+x}, \quad h_x-x+1\leq i \leq n-r-x.
\end{equation}
\item If (\ref{eqdc}) holds,  then 
\begin{equation}\label{eqdcgeq}
d_i\geq c_{i+x}, \quad 1\leq i \leq n-r-x.
\end{equation}
\end{enumerate}

\end{lemma}

\begin{rem}\label{remlemmahx}
If $x=0$,  then $h_0=0$ and (\ref{eqdc}) is equivalent to $\bd=\bc$.
\end{rem}

\noindent{\bf Proof of Lemma \ref{lemmahx}.}
\begin{enumerate}
\item
If $\bc\prec'(\bd, \ba)$, then
(\ref{eqdcgeq}) and
$$
\sum_{i=1}^{h_x}c_i- \sum_{i=1}^{h_x-x}d_i\leq \sum_{i=1}^{x}a_i=\sum_{i=1}^{n-r}c_i- \sum_{i=1}^{n-r-x}d_i
$$
hold. Thus, $$0\leq \sum_{i=h_x+1}^{n-r}c_i- \sum_{i=h_x-x+1}^{n-r-x}d_i=\sum_{i=h_x-x+1}^{n-r-x}(c_{i+x}- d_i).$$
  From (\ref{eqdcgeq}) we obtain (\ref{eqdc}).
\item Assume that (\ref{eqdc}) holds. From the definition of $h_x$, 
   $d_{i-x+1}\geq c_i$   for $x\leq i \leq  h_x-1$. Therefore,
   $d_i\geq c_{i+x-1}\geq c_{i+x}$ for $1\leq i \leq  h_x-x$, and the result follows.
\end{enumerate}
\hfill $\Box$

\subsection{Previous results}\label{secprevious}

Along this paper,  we repeatedly introduce different  sequences of integers $\ba=(a_1,\dots, a_x)$ and $\bb=(b_1,\dots,b_{z-x})$. By Lemma \ref{lemmadec}, $\ba$ and $\bb$ will be well defined and $\bb$ will be a partition, that is, $a_1\geq\dots\geq a_x$ and $b_1\geq\dots\geq b_{z-x}\geq 0$. 

\begin{lemma}\label{lemmadec}
 Let $z,x$ be integers such that $0\leq x\leq z$, and let  
$ \phi_1(s)\mid\cdots\mid \phi_{  r}(s)$ and
$\gamma_1(s)\mid\cdots\mid\gamma_{  r+x}(s)$ be two polynomial chains. Then
 $$
\begin{array}{rl}
     &\deg(\gamma_{x-j+1})+\sum_{i=1}^{  r}\deg(\lcm(    \phi_{i},
  \gamma_{i+x-j+1}))-\sum_{i=1}^{    r}\deg(\lcm(    \phi_{i},   \gamma_{i+x-j}))
\\ \geq &\deg(\gamma_{x-j})+\sum_{i=1}^{    r}\deg(\lcm(    \phi_{i},   \gamma_{i+x-j}))-\sum_{i=1}^{    r}\deg(\lcm(    \phi_{i},
  \gamma_{i+x-j-1})),\\&\hfill 1\leq j\leq x-1,
\end{array}
$$
and 
$$
\begin{array}{rl}
&-\deg(\gamma_{x+j})+\sum_{i=1}^{    r-j+1}\deg(\lcm(   \phi_{i},   \gamma_{i+x+j-1}))-\sum_{i=1}^{ 
  r-j}\deg(\lcm(   \phi_{i},   \gamma_{i+x+j}))\\\geq&-\deg(\gamma_{x+j+1})+
\sum_{i=1}^{    r-j}\deg(\lcm(   \phi_{i},   \gamma_{i+x+j}))-\sum_{i=1}^{r-j-1}\deg(\lcm(   \phi_{i},   \gamma_{i+x+j+1}))\\ &\hfill \quad 1\leq j\leq z-x-1.
\end{array}
$$
As a consequence, for nonnegative integers 
$p_1\leq\dots\leq p_r$ and $q_1\leq\dots\leq q_{r+x}$,
$$
\begin{array}{rl}
     &q_{x-j+1}+\sum_{i=1}^r \max\{p_i,q_{i+x-j+1}\}-\sum_{i=1}^r \max\{p_i,q_{i+x-j}\}
     \\\geq &q_{x-j}+\sum_{i=1}^r \max\{p_i,q_{i+x-j}\}-\sum_{i=1}^r \max\{p_i,q_{i+x-j-1}\},\quad  1\leq j\leq x-1,
\end{array}
$$
and
$$
\begin{array}{rl}
     &-q_{x+j}+\sum_{i=1}^{r-j+1} \max\{p_i,q_{i+x+j-1}\}-\sum_{i=1}^{r-j} \max\{p_i,q_{i+x+j}\}
     \\\geq &-q_{x+j+1}+\sum_{i=1}^{r-j} \max\{p_i,q_{i+x+j}\}-\sum_{i=1}^{r-j-1} \max\{p_i,q_{i+x+j+1}\},\quad  1\leq j\leq z-x-1.
\end{array}
$$
\end{lemma}

 {\noindent\bf Sketch of the proof.}
The first inequality follows from expression (4) of \cite[Lemma 2]{DoStEJC10} applied to
$ \phi_1(s)\mid\cdots\mid \phi_{  r}(s)$ and
$\gamma_1(s)\mid\cdots\mid\gamma_{  r+x}(s)$ (and replacing $j$ by $x-j$). In order to obtain the second inequality, define 
$$
\widehat{\phi}_i(s)=\phi_{i-z}(s), \quad 1\leq i \leq r+z,
$$
apply \cite[Lemma 2]{DoStEJC10} to
$\gamma_1(s)\mid\cdots\mid\gamma_{  r+x}(s)$ and
$\widehat\phi_1(s)\mid\cdots\mid\widehat \phi_{  r+z}(s)$,  afterwards replace $j$ by $z-x-j$,  and then $i$ by $i+z$. 
\hfill $\Box$

\medskip
The following lemma 
relates the 
complete structural data of a rational matrix with the complete structural data of certain polynomial matrices. 
 \begin{lemma}[{\rm\cite[Lemma 2.3]{AmBaMaRo25}}]\label{lem_polrat}
Let $R(s)$ be a rational matrix and let $p(s)$ be  a  monic polynomial multiple	of the least common denominator of the entries in $R(s)$. Then, $p(s)R(s)$ is a polynomial matrix of the same rank as $R(s)$ and:
\begin{itemize}
	\item[(i)] The quotients  $\frac{\epsilon_1(s)}{\psi_1(s)},\dots,\frac{\epsilon_r(s)}{\psi_r(s)}$ are the invariant rational functions of $R(s)$ if and only if the polynomials $\frac{p(s)\epsilon_1(s)}{\psi_1(s)},\dots,\frac{p(s)\epsilon_r(s)}{\psi_r(s)}$ are the invariant factors of $p(s)R(s)$. 
	\item[(ii)]
 The integers $\tilde q_1,\dots,\tilde q_r$ are the invariant orders at $\infty$ of $R(s)$ if and only if 
 $\tilde q_1-\deg(p(s)), \dots,\tilde q_r-\deg(p(s))$ are the invariant orders at $\infty$ of  $p(s)R(s)$. 
	\item[(iii)] $\mathcal{N}_r(R(s))=\mathcal{N}_r(p(s)R(s))$, $\mathcal{N}_{\ell}(R(s))=\mathcal{N}_{\ell}(p(s)R(s))$ and, therefore, the minimal indices of $p(s)R(s)$ and of $R(s)$ are the same.
\end{itemize}
\end{lemma}

In Theorem \ref{theoexistencerat} we  find  necessary and sufficient conditions for the existence of a rational matrix  with prescribed complete structural data. This result was proved in \cite[Theorem 4.1]{AnDoHoMa19} for infinite fields and in \cite[Theorem 2.4]{AmBaMaRo25} for arbitrary fields.

\begin{theorem}[{\rm\cite[Theorem 2.4]{AmBaMaRo25}}]\label{theoexistencerat}
Let $m, n, r\leq \min\{m,n\}$  be  positive integers.
Let
$\epsilon_1(s)\mid\dots \mid  \epsilon_r(s)$ and 
 $\psi_r(s)\mid\dots\mid \psi_1(s)$ be monic polynomials
such that  $\frac{\epsilon_1(s)}{\psi_1(s)},\dots,\frac{\epsilon_r(s)}{\psi_r(s)}$ are irreducible rational functions.
Let $\tilde q_1\leq \dots \leq \tilde q_r$ be  integers and $(d_1, \ldots, d_{n-r})$, $(v_1, \ldots, v_{m-r})$
 partitions. 
 There exists a rational matrix $R(s)\in \FF(s)^{m\times n}$, $\rank(R(s))= r$,  with  $\frac{\epsilon_1(s)}{\psi_1(s)},\dots,\frac{\epsilon_r(s)}{\psi_r(s)}$ as invariant rational functions,   $\tilde q_1, \dots, \tilde q_r$  as invariant orders at $\infty$, and $d_1, \dots, d_{n-r}$ and $v_1, \dots, v_{m-r}$ as column and row minimal indices, respectively, if and only if
\begin{equation}\label{eqISTrat}
\sum_{i=1}^{n-r}d_i+\sum_{i=1}^{m-r}v_i+\sum_{i=1}^{r}\tilde q_i+\sum_{i=1}^{r}\deg(\epsilon_i)-\sum_{i=1}^{r}\deg(\psi_i)=0.
\end{equation}
\end{theorem}

From now on we use the following notation:
given $\varphi(s),\eta(s),\psi(s),\epsilon(s)\in\efe[s]$ such that  $\gcd(\varphi,\eta)=1$ and $\gcd(\psi,\epsilon)=1$, and $p, q$ integers, we denote
$$
\Delta\left(\frac{\eta}{\varphi}, \frac{\epsilon}{\psi}, p, q\right)=
\deg(\lcm(\eta,\epsilon))- \deg(\gcd(\varphi,\psi)) +\max\{p,q\},
$$
$$\Delta\left(\frac{\eta}{\varphi}, \frac{\epsilon}{\psi}\right)=
\deg(\lcm(\eta,\epsilon))- \deg(\gcd(\varphi,\psi)),$$
$$\Delta\left(\frac{\eta}{\varphi}, p\right)=\deg(\eta)- \deg(\varphi)+p,$$

$$\Delta\left(\frac{\eta}{\varphi}\right)=
\deg(\eta)- \deg(\varphi).$$

The next theorem contains a solution to the row completion problem for rational matrices when  the complete structural data are precribed. 

\begin{theorem}[{\rm\cite[Theorem 3.11]{AmBaMaRo25}}] {\rm (Prescription of the complete structural data for rational matrices)}\label{prescr4rat}
  Let $R(s)\in\efe(s)^{m\times n}$ be a rational matrix, $\rank (R(s))=r$. 
 Let $\frac{\eta_1(s)}{\varphi_1(s)},\dots, \frac{\eta_r(s)}{\varphi_r(s)}$ be its invariant rational functions,  
$\tilde p_1, \dots, \tilde p_r$ its invariant orders at $\infty$, 
	$\bc=(c_1,  \dots,  c_{n-r})$  its
	column minimal indices, and $\bu=(u_1, \dots,  u_{m-r})$  its row minimal indices.
	
	Let $z,x$ be integers such that $0\leq x\leq \min\{z, n-r\}$ and let $\epsilon_1(s)\mid \dots \mid\epsilon_{r+x}(s)$ and 
        $\psi_{r+x}(s)\mid \dots \mid\psi_{1}(s)$ be monic polynomials such that $\frac{\epsilon_i(s)}{\psi_i(s)}$ are irreducible rational functions, $1\leq i \leq r+x$.
        Let $\tilde q_{1}\leq \dots \leq \tilde q_{r+x}$ be integers and $\bd=(d_1,  \dots,  d_{n-r-x})$  and $\bv=(v_1, \dots,  v_{m+z-r-x})$  be two partitions.
	There exists a rational matrix  $\widetilde W(s)\in \efe(s)^{z\times n}$ such that $\rank \left(\begin{bmatrix}R(s)\\\widetilde W(s)\end{bmatrix}\right)=r+x$ and
        $\begin{bmatrix}R(s)\\\widetilde W(s)\end{bmatrix}$ has $\frac{\epsilon_1(s)}{\psi_1(s)}, \dots, \frac{\epsilon_{r+x}(s)}{\psi_{r+x}(s)}$ as invariant rational functions, 
 $\tilde q_1, \dots,\tilde q_{r+x}$ as invariant orders at $\infty$,  $d_1,  \dots,  d_{n-r-x}$  as column minimal indices and $v_1, \dots,  v_{m+z-r-x}$ as row minimal indices  
        if and only if 
        \begin{equation}\label{eqinterratnum}
		\epsilon_i(s)\mid \eta_i(s)\mid \epsilon_{i+z}(s),\quad 1\leq i \leq r,
	\end{equation}
        \begin{equation}\label{eqinterratden}
		\psi_{i+z}(s)\mid \varphi_i(s)\mid \psi_{i}(s),\quad 1\leq i \leq r,
	\end{equation}
        \begin{equation}\label{eqinterratioi}
\tilde q_i\leq \tilde p_i\leq \tilde q_{i+z}, \quad 1\leq i\leq r,
        \end{equation}
	\begin{equation}\label{eqcmimajrat}  \bc \prec'  (\bd , \tilde \ba),\end{equation}
	\begin{equation}\label{eqrmimajrat}\bv \prec'  (\bu , \tilde \bb),\end{equation}
\begin{equation}\label{eqdegsumrat}
\begin{aligned}
  \sum_{i=1}^{r}
  \Delta\left(\frac{\eta_i}{\varphi_i}, \frac{\epsilon_{i+x}}{\psi_{i+x}}, \tilde p_i,\tilde q_{i+x}\right)-\sum_{i=1}^{r}\Delta\left(\frac{\epsilon_{i+x}}{\psi_{i+x}}, \tilde q_{i+x}\right)
  \leq \sum_{i=1}^{m+z-r-x}v_i-\sum_{i=1}^{m-r}u_i
  ,\\
  \mbox{ with equality when $x=0$,}
\end{aligned}
 \end{equation}
	where $\tilde \ba = (\tilde a_1, \dots, \tilde a_x )$ and $\tilde \bb = (\tilde b_1, \dots,  \tilde b_{z-x} )$ are defined as
 \begin{equation}\label{eqdeftildea}
		\begin{array}{rl}
		\sum_{i=1}^{j}	\tilde a_i=&
			\sum_{i=1}^{m+z-r-x}v_i-\sum_{i=1}^{m-r} u_i+\sum_{i=1}^{r+j}\Delta\left(\frac{\epsilon_{i+x-j}}{\psi_{i+x-j}}, \tilde q_{i+x-j}\right)\\&
-\sum_{i=1}^{  r} \Delta\left(\frac{\eta_i}{\varphi_i}, \frac{\epsilon_{i+x-j}}{\psi_{i+x-j}}, \tilde p_{i}, \tilde q_{i+x-j}\right),\quad
                        1\leq j \leq x,
		\end{array}
\end{equation}
\begin{equation}\label{eqdeftildeb}
	\begin{array}{rl}
 \sum_{i=1}^{j}\tilde b_i=&
			\sum_{i=1}^{m+z-r-x}  v_i-\sum_{i=1}^{m-r}  u_i+\sum_{i=1}^{  r-j}\Delta\left(\frac{\epsilon_{i+x+j}}{\psi_{i+x+j}},\tilde q_{i+x+j}\right)\\&
-\sum_{i=1}^{  r-j}
\Delta\left(\frac{\eta_i}{\varphi_i}, \frac{\epsilon_{i+x+j}}{\psi_{i+x+j}}, \tilde p_i,  \tilde q_{i+x+j}\right),\quad
                        1\leq j \leq z-x,\\
		\end{array}
	\end{equation}  
\end{theorem}

\begin{rem}\label{remdeceta}\ 
    \begin{enumerate}
     \item\label{remdeceta.eta}
     Let $\eta = \#\{i : u_i > 0\}$ and $\bar \eta = \#\{i : v_i > 0\}$.
     Theorem 3.8  in \cite{AmBaMaRo25} also contains the condition
\begin{equation}\label{eqetapol}
\bar \eta \geq  \eta.
\end{equation}
But in \cite[Lemma 1]{DoSt25} it is shown that this condition is redundant ((\ref{eqrmimajrat}) implies (\ref{eqetapol})). 

\item \label{remdeceta.remratcdioi}
By  Theorem \ref{theoexistencerat}  (see also \cite[Remark 4.3]{AmBaMaRo24}),
if (\ref{eqinterratnum})--(\ref{eqdegsumrat})
 hold, then
 \begin{equation}\label{eqdegsumratcdioi}
\begin{aligned}
  \sum_{i=1}^{r}\Delta\left(\frac{\eta_i}{\varphi_i}, \frac{\epsilon_{i+x}}{\psi_{i+x}}, \tilde p_i,\tilde q_{i+x}\right) 
  \leq \sum_{i=1}^{n-r}c_i-\sum_{i=1}^{n-r-x}d_i
  +\sum_{i=1}^{r}\Delta\left(\frac{\eta_i}{\varphi_i}, \tilde p_i\right) 
  -  \sum_{i=1}^{x}\Delta\left(\frac{\epsilon_{i}}{\psi_{i}},\tilde q_{i}\right),\\
\mbox{ with equality when $x=z$,} 
  \end{aligned}
 \end{equation}
and (\ref{eqcmimajrat}) and  (\ref{eqrmimajrat}) hold for
   $\tilde \ba=(\tilde a_1, \dots, \tilde a_x)$ and $\tilde \bb=(\tilde b_1, \dots, \tilde b_{z-x})$ defined as
   \begin{equation}\label{eqdeftildeacdioi}
                \begin{array}{rl}
			 \sum_{i=1}^{j}\tilde a_i=& \sum_{i=1}^{n-r}c_i-\sum_{i=1}^{n-r-x}d_i+
                        \sum_{i=1}^{r}\Delta\left(\frac{\eta_i}{\varphi_i}, \tilde p_i\right) 
  -  \sum_{i=1}^{x-j}\Delta\left(\frac{\epsilon_{i}}{\psi_{i}},\tilde q_{i}\right)
  \\&-\sum_{i=1}^{r}\Delta\left(\frac{\eta_i}{\varphi_i}, \frac{\epsilon_{i+x-j}}{\psi_{i+x-j}}, \tilde p_i,\tilde q_{i+x-j}\right),
               \quad 1\leq j \leq x,
		\end{array}
\end{equation}
\begin{equation}\label{eqdeftildebcdbioi}
	\begin{array}{rl}
  \sum_{i=1}^{j}\tilde b_i=&
\sum_{i=1}^{n-r}c_i-\sum_{i=1}^{n-r-x}d_i
+
                        \sum_{i=1}^{r}\Delta\left(\frac{\eta_i}{\varphi_i}, \tilde p_i\right)-\sum_{i=1}^{x+\min\{j,r\}}\Delta\left(\frac{\epsilon_{i}}{\psi_{i}},\tilde q_{i}\right) \\
   &-\sum_{i=1}^{r-j}\Delta\left(\frac{\eta_i}{\varphi_i}, \frac{\epsilon_{i+x+j}}{\psi_{i+x+j}}, \tilde p_i,\tilde q_{i+x+j}\right),
                        \quad
                        1\leq j \leq z-x.
		\end{array}
	\end{equation}
Conversely, (\ref{eqinterratnum})--(\ref{eqrmimajrat}) and 
(\ref{eqdegsumratcdioi}) 
with $\tilde \ba$ and $\tilde \bb$   defined as in 
(\ref{eqdeftildeacdioi}) and (\ref{eqdeftildebcdbioi}), respectively, imply
 (\ref{eqinterratnum})--(\ref{eqdegsumrat}) 
 with $\tilde \ba$ and $\tilde \bb$   defined as in 
(\ref{eqdeftildea}) and (\ref{eqdeftildeb}), respectively.

    \end{enumerate}

\end{rem}

The next result is a particular case  of Theorem \ref{prescr4rat} for polynomial matrices. In fact, if $R(s)=P(s)$ is a polynomial matrix,  i.e., $\varphi_1(s)=1$, and we prescribe $\psi_1(s)=1$, then the completed matrix must be polynomial. 

\begin{corollary}[{\rm\cite[Theorem 3.8]{AmBaMaRo25}}] {\rm (Prescription of the complete structural data for  polynomial matrices)}
  \label{theoprescr4ioi}
  Let $P(s)\in\efe[s]^{m\times n}$ be a polynomial matrix, $\rank (P(s))=r$.
	Let $\alpha_1(s)\mid \cdots\mid\alpha_{r}(s)$ be its 
invariant factors, $p_1, \dots, p_r$ its invariant orders at $\infty$,
	$\bc=(c_1,  \dots,  c_{n-r})$  its
	column minimal indices, and $\bu=(u_1, \dots,  u_{m-r})$  its row minimal indices.
	 
  Let $z,x$ be integers such that $0\leq x\leq \min\{z, n-r\}$. 
	 	Let $\beta_1(s)\mid \cdots\mid\beta_{r+x}(s)$ be monic polynomials, $q_1\leq\dots\leq q_{r+x}$ be integers, and 
	$\bd=(d_1,  \dots,  d_{n-r-x})$  and $\bv=(v_1, \dots,  v_{m+z-r-x})$ be two partitions.
	There exists a polynomial matrix $W(s)\in \efe[s]^{z\times n}$ such that  	
 $\rank \left(\begin{bmatrix}P(s)\\W(s)\end{bmatrix}\right) =r+x$ and 
	$\begin{bmatrix}P(s)\\W(s)\end{bmatrix}$ has $\beta_1(s), \dots, \beta_{r+x}(s)$ as invariant factors,
	$q_1, \dots, q_{r+x}$ as invariant orders at $\infty$,
 $d_1, \dots, d_{n-r-x}$ as column minimal indices
	and 
	$v_1, \dots, v_{m+z-r-x}$ as row minimal indices  
	if and only if
\begin{equation}\label{eqinterif}
        \beta_i(s)\mid \alpha_{i}(s)\mid \beta_{i+z}(s),\quad 1\leq i \leq r,
 \end{equation}
\begin{equation}\label{eqinterpolioi}
q_i\leq p_i\leq q_{i+z}, \quad 1\leq i\leq r,
        \end{equation}
	\begin{equation}\label{eqcmimajpol}  \bc \prec'  (\bd , \ba),\end{equation}
	\begin{equation}\label{eqrmimajpol}\bv \prec'  (\bu , \bb),\end{equation}
	\begin{equation}\label{eqdegsumpolioi}
   \begin{aligned}
  \sum_{i=1}^{r}\deg(\lcm(\alpha_i,\beta_{i+x}))+
  \sum_{i=1}^{r}\max\{p_i,q_{i+x}\}\\
  \leq \sum_{i=1}^{m+z-r-x}v_i-\sum_{i=1}^{m-r}u_i
  +\sum_{i=1}^{r}\deg(\beta_{i+x}) +\sum_{i=1}^{r}q_{i+x},\\
 \mbox{ with equality when $x=0$,}
\end{aligned}
 \end{equation}
	where $ \ba = (a_1, \dots, a_x )$ and $\bb = ( b_1, \dots,  b_{z-x} )$ are defined as
	\begin{equation}\label{eqdefaioi}
 \begin{array}{rl}
			  \sum_{i=1}^{j}a_i=&
			\sum_{i=1}^{m+z-r-x}v_i-\sum_{i=1}^{m-r} u_i+
   \sum_{i=1}^{ r+j}\deg(  \beta_{i+x-j})+\sum_{i=1}^{ r+j}q_{i+x-j}
                        \\&
   -\sum_{i=1}^{ r}\deg(\lcm( \alpha_{i},  \beta_{i+x-j}))
                        -\sum_{i=1}^{r}\max\{ p_i, q_{i+x-j}\},\\&\hfill
			1\leq j \leq x,
		\end{array}
\end{equation}
\begin{equation}\label{eqdefbioi}
		\begin{array}{rl}
\sum_{i=1}^{j} b_i=&
\sum_{i=1}^{m+z-r-x}v_i-\sum_{i=1}^{m-r} u_i
                        +\sum_{i=1}^{r-j}\deg(\beta_{i+x+j})+\sum_{i=1}^{r-j}
                        q_{i+x+j}\\&
			-\sum_{i=1}^{r-j}\deg(\lcm(\alpha_i, \beta_{i+x+j}))-
\sum_{i=1}^{r-j}\max\{ p_i,  q_{i+x+j}\},\\& \hfill
  1 \leq j \leq z-x.
		\end{array}
  \end{equation}

\end{corollary}

\section{Row (column) completion for polynomial and rational matrices: remaining cases}\label{secremaining}

The aim of this section is to give a solution to Problem \ref{problemrat} for the cases not solved in \cite{AmBaMaRo25}. We present the solution in different subsections depending on the invariant(s) prescribed.
Several proofs are analogous to some proofs in \cite{AmBaMaRo24_2}. For the reader's convenience we have  included them in  \ref{secappendix}.

\subsection{Prescription of the infinite and singular structures}\label{subsecinfsing}

In this subsection we prescribe the infinite and singular structures for both  polynomial and rational matrices. 

The proof of Theorem \ref{theoprescrioirmicmipol}
 is analogous to that of \cite[Theorem 3.1]{AmBaMaRo24_2}. It can be found in \ref{secappendix}.

\begin{theorem} {\rm (Prescription of the infinite and singular structures for polynomial matrices)} \label{theoprescrioirmicmipol}
  Let $\efe$ be an algebraically closed field.
  Let $P(s)\in\efe[s]^{m\times n}$ be a polynomial matrix, $\rank (P(s))=r$.
Let $\alpha_1(s)\mid\cdots\mid\alpha_r(s)$ be the invariant factors,
$p_1, \dots,  p_r$ the invariant orders at $\infty$,
$\bc=(c_1,  \dots,  c_{n-r})$  the
column minimal indices, and $ \bu=(u_1, \dots,  u_{m-r})$  the row minimal
indices of $P(s)$.

Let $z$ and $x$ be integers such that $0\leq x\leq \min\{z, n-r\}$ and 
let  $q_1\leq\dots \leq q_{r+x}$ be integers,  and 
 $\bd=(d_1, \dots, d_{n-r-x})$ and $\bv=(v_1, \dots, v_{m+z-r-x})$ be partitions.
Let
\begin{equation}\label{eqbigAioicmirmi}
A=\sum_{i=1}^{r+x}q_i-\sum_{i=1}^{r}p_i+\sum_{i=1}^{n-r-x}d_i-\sum_{i=1}^{n-r}c_i+\sum_{i=1}^{m+z-r-x}v_i-\sum_{i=1}^{m-r}u_i.
\end{equation}
There exists $W(s)\in \efe[s]^{z\times n}$ such that  $\rank \left(\begin{bmatrix}P(s)\\W(s)\end{bmatrix}\right)=r+x$
and $\begin{bmatrix}P(s)\\W(s)\end{bmatrix}$  has 
 $q_1,\dots, q_{r+x} $ as  invariant orders at $\infty$, 
 $\bd=(d_1,\dots, d_{n-r-x})$ as  column minimal indices and
$\bv=(v_1,\dots, v_{m+z-r-x})$ as row minimal indices
if and only if 
 (\ref{eqinterpolioi}), 
\begin{equation}\label{eqAleq}
A\leq \sum_{i=r+x-z+1}^{r}\deg(\alpha_i),
\end{equation}
\begin{equation}\label{eqvusAioi}
 \sum_{i=1}^{m+z-r-x}v_i-\sum_{i=1}^
 {m-r}u_i\geq \max\{0, A\}+\sum_{i=1}^{r}\max\{p_{i}, q_{i+x}\}-\sum_{i=1}^{r}q_{i+x},
\end{equation}
\begin{equation}\label{eqcmimajpolhatsAioi}  \bc \prec'  (\bd , \hat \ba),\end{equation}
\begin{equation}\label{eqrmimajpolhatsAioi} \bv \prec'  (\bu , \hat \bb),\end{equation}
where $\hat \ba = (\hat a_1, \dots, \hat a_x )$ and $\hat \bb = (\hat b_1, \dots,  \hat b_{z-x})$ 
are  defined as
\begin{equation}\label{eqaioicmirmi}
\sum_{i=1}^{j}\hat a_i=
  \sum_{i=1}^{m+z-r-x}v_i-\sum_{i=1}^{m-r}u_i+\sum_{i=1}^{r+j}q_{i+x-j}-\sum_{i=1}^{r}\max\{p_{i}, q_{i+x-j}\}-A,
  \quad 1\leq j \leq x,
\end{equation}
\begin{equation}\label{eqbioicmirmi}
\begin{array}{rl}
 \sum_{i=1}^{j} \hat b_i=&
			\sum_{i=1}^{m+z-r-x}  v_i-\sum_{i=1}^{m-r}  u_i+
   \min\{0,\sum_{i=r-j+1}^{r}\deg(\alpha_{i})-A\}                       
\\&+\sum_{i=1}^{  r-j}q_{i+x+j}-
\sum_{i=1}^{  r-j}\max\{p_i,  q_{i+x+j}\}, \quad  1\leq j \leq z-x.
		\end{array}
\end{equation}

\end{theorem}

\begin{rem}\label{remacdecr} Analogous to \cite[Remark 3.2]{AmBaMaRo24_2}. 
  \begin{enumerate}
    \item\label{remitac}
The necessity of the conditions,  and the sufficiency when $A\leq0$ hold for arbitrary fields.   
\item\label{remitemhatbdecr}
        Let us see that  if (\ref{eqAleq}) and (\ref{eqvusAioi}) hold, then  $\hat b_1\geq \dots \geq \hat b_{z-x}\geq 0$.
        By (\ref{eqAleq}), $z-x\in \{k\geq 0: A\leq \sum_{i=r-k+1}^{r}\deg(\alpha_{i}) \}$. Let $g=\min\{k\geq 0: A\leq \sum_{i=r-k+1}^{r}\deg(\alpha_{i}) \}$.
        Then $0\leq g \leq z-x$. Observe that  $g=0$ if and only if $A\leq 0$. If $g\geq 1$, then
        \begin{equation}\label{eqtarte}
          \sum_{i=r-g+2}^{r}\deg(\alpha_{i}) <A\leq \sum_{i=r-g+1}^{r}\deg(\alpha_{i}).\end{equation}
        If $g=0$ or $g=1$, then $A\leq \deg(\alpha_r)\leq  \sum_{i=r-j+1}^{r}\deg(\alpha_{i}) $ for $1\leq j \leq z-x$; hence
$$
 \sum_{i=1}^{j} \hat b_i=
			\sum_{i=1}^{m+z-r-x}  v_i-\sum_{i=1}^{m-r}  u_i                      +\sum_{i=1}^{  r-j}q_{i+x+j}-
\sum_{i=1}^{  r-j}\max\{p_i,  q_{i+x+j}\}, \quad 1\leq j \leq z-x,
$$

If $g\geq 2$, then from (\ref{eqtarte}) we get
        $$
        \begin{array}{l}
        \sum_{i=r-j+1}^{r}\deg(\alpha_{i})\leq \sum_{i=r-g+2}^{r}\deg(\alpha_{i}) <A, \quad 1\leq j \leq g-1,\\ \\
        A\leq \sum_{i=r-g+1}^{r}\deg(\alpha_{i})\leq \sum_{i=r-j+2}^{r}\deg(\alpha_{i})\leq \sum_{i=r-j+1}^{r}\deg(\alpha_{i}),\\\hfill g+1\leq j \leq z-x.
        \end{array}
        $$
        Thus,  
$$      
\begin{array}{rl}
 \sum_{i=1}^{j} \hat b_i=&
			\sum_{i=1}^{m+z-r-x}  v_i-\sum_{i=1}^{m-r}  u_i+
   \sum_{i=r-j+1}^{r}\deg(\alpha_{i})-A                     
\\&+\sum_{i=1}^{  r-j}q_{i+x+j}-
\sum_{i=1}^{  r-j}\max\{p_i,  q_{i+x+j}\}, \quad  1\leq j \leq g-1,\\ \\
\sum_{i=1}^{j} \hat b_i=&
			\sum_{i=1}^{m+z-r-x}  v_i-\sum_{i=1}^{m-r}  u_i
+\sum_{i=1}^{  r-j}q_{i+x+j}\\&
-
\sum_{i=1}^{  r-j}\max\{p_i,  q_{i+x+j}\}, \quad  g\leq j \leq z-x.\\
		\end{array}
$$
        From (\ref{eqvusAioi}) and (\ref{eqtarte}), we  have
        $\sum_{i=1}^{m+z-r-x}v_i-\sum_{i=1}^{m-r}u_i-A\geq \sum_{i=1}^{r}\max\{p_{i}, q_{i+x}\}-\sum_{i=1}^{r}q_{{i+x}}$
         and
        $0<A-\sum_{i=r-g+2}^{r}\deg(\alpha_{i})\leq  \deg(\alpha_{r-g+1})$, respectively.
       Thus, by Lemma \ref{lemmadec},
    we  obtain   $\hat b_1\geq \dots \geq \hat b_{z-x}\geq 0$.
  \end{enumerate}
\end{rem}

In the following example we show that, in general, conditions (\ref{eqinterpolioi}) and (\ref{eqAleq})-(\ref{eqrmimajpolhatsAioi}) are not sufficient if the field is not algebraically closed and $A>0$.
Observe also that, to achieve the desired invariants, the degree of the completed matrix must be necessarily greater than the degree of the prescribed submatrix.

\begin{example}
  \label{exalgclosed}

Let $P(s)=
  \begin{bmatrix}0&s&1\\s^2+1&0&0
  \end{bmatrix}
  \in \FF[s]^{2\times3}
  $.
  We have
  $$
  r=2, \quad \alpha_1(s)=1, \quad \alpha_2(s)=s^2+1,  \quad p_1=-2, \quad p_2=-1, \quad \bc=(1), \quad\bu=\emptyset.$$ 
 Let $x=0$, $z=1$.
 We prescribe
$$
q_1, q_2 \text{ integers such that } q_1\leq q_2, \quad \bd=(1), \quad \bv=(5).$$ 
 Then, $A=q_1+q_2+8$ and it is easy to see that (\ref{eqinterpolioi}) and (\ref{eqAleq})--(\ref{eqrmimajpolhatsAioi}) hold if and only if $ q_1 \leq -2\leq  q_2\leq -1$ and  $-8\leq q_1+q_2\leq -6$.
Therefore, we must prescribe $q_1\leq -4$, hence   
the degree of the completed polynomial matrix is greater than or equal to 4, that is, it is greater than the degree of $P(s)$. 
\begin{itemize}
\item If $q_1=-6$, $q_2=-1$, then  
$A=1>0$. If $\FF$ is an algebraically closed field, then by Theorem  \ref{theoprescrioirmicmipol} there exists $W(s)\in \efe[s]^{1\times 3}$ such that   $\begin{bmatrix}P(s)\\W(s)\end{bmatrix}$
has the prescribed invariants.
Observe that if 
  $\beta_1(s) \mid \beta_2(s)$ are its invariant factors then, by Theorem \ref{theoexistencerat},
  $\deg(\beta_1)+ \deg(\beta_2)=1$, and by Corollary \ref{theoprescr4ioi}, $\beta_1(s) \mid \alpha_1(s)=1$ and
  $\beta_2(s) \mid \alpha_2(s)=s^2+1$, i.e., there exists a
   polynomial $\beta_2(s)$ such that $\beta_2(s) \mid
  s^2+1$ and $\deg(\beta_2)=1$.
  If $\efe=\RR$, there is no such polynomial.   If  $\FF=\CC$, a possible completion is
$$
  \begin{bmatrix}P(s)\\W(s)\end{bmatrix}=
  \begin{bmatrix}0&s&1\\s^2+1&0&0\\s^5(s+i)&0&0
  \end{bmatrix}
  \in \CC[s]^{3\times 3}.$$
In particular, we see that if $\FF$ is not algebraically closed and $A>0$, conditions 
(\ref{eqinterpolioi}) and (\ref{eqAleq})--(\ref{eqrmimajpolhatsAioi}) are no sufficient to guarantee the existence of the desired completion.
\item If $q_1=-6$, $q_2=-2$, then $A=0$. The matrix
  $$\begin{bmatrix}P(s)\\W(s)\end{bmatrix}=
 \begin{bmatrix}0&s&1\\s^2+1&0&0\\
 0&s^6&s^5
  \end{bmatrix}
  \in \FF[s]^{3\times 3} 
  $$ has the prescribed invariants  for an arbitrary field $\FF$.

\end{itemize}
\end{example}

In the following theorem, as in Theorem \ref{theoprescrioirmicmipol}, we solve Problem \ref{problemrat} when the infinite structure and the column and row minimal indices are prescribed, but for rational matrices. In this case the field is not required to be algebraically closed.

\begin{theorem} {\rm (Prescription of the infinite and singular structures for rational matrices)}
  \label{theoprescrioirmicmirat}
  Let $R(s)\in\efe(s)^{m\times n}$ be a rational matrix, $\rank (R(s))=r$. 
 Let 
$\tilde p_1, \dots, \tilde p_r$ be its invariant orders at $\infty$, 
	$\bc=(c_1,  \dots,  c_{n-r})$  its
	column minimal indices, and $\bu=(u_1, \dots,  u_{m-r})$  its row minimal
	indices.
	
	Let $z,x$ be integers such that $0\leq x\leq \min\{z, n-r\}$ and
        let $\tilde q_{1}\leq \dots \leq \tilde q_{r+x}$ be integers, and $\bd=(d_1,  \dots,  d_{n-r-x})$  and $\bv=(v_1, \dots,  v_{m+z-r-x})$  be two partitions.
       Let $h_x=\min\{i: d_{i-x+1}<c_i\}$.    
      There exists a rational matrix $\widetilde W(s)\in \efe(s)^{z\times n}$ such that $\rank \left(\begin{bmatrix}R(s)\\\widetilde W(s)\end{bmatrix}\right)=r+x$, and
        $\begin{bmatrix}R(s)\\\widetilde W(s)\end{bmatrix}$ has 
 $\tilde q_1, \dots,\tilde q_{r+x}$ as invariant orders at $\infty$,  $d_1,  \dots,  d_{n-r-x}$  as column minimal indices and $v_1, \dots,  v_{m+z-r-x}$ as row minimal indices  
        if and only if 
        (\ref{eqdc}), (\ref{eqinterratioi}), 
 \begin{equation}\label{eqvusAioirat}
 \sum_{i=1}^{m+z-r-x}v_i-\sum_{i=1}^{m-r}u_i\geq \sum_{i=1}^{r}\max\{\tilde p_{i}, \tilde q_{i+x}\}-\sum_{i=1}^{r}\tilde q_{i+x},
\end{equation}
\begin{equation}\label{eqrmimajrathatsAioi} \bv \prec'  (\bu , \hat{\tilde \bb}),\end{equation}
where 
$\hat{\tilde \bb}=(\hat{\tilde b}_1, \dots, \hat{\tilde b}_{z-x})$ is defined as
\begin{equation}\label{eqbioicmirmirat}
\begin{array}{rl}
  \sum_{i=1}^j\hat {\tilde b}_i=&
			\sum_{i=1}^{m+z-r-x}  v_i-\sum_{i=1}^{m-r}  u_i                    
+\sum_{i=1}^{  r-j}\tilde q_{i+x+j}\\&-
\sum_{i=1}^{  r-j}\max\{\tilde p_i, \tilde q_{i+x+j}\},
 \quad  1\leq j \leq z-x.
		\end{array}
\end{equation}
\end{theorem}

{\noindent\bf Proof.}
Assume that there exists $\widetilde W(s)\in \efe(s)^{z\times n}$ such that  $\rank \left(\begin{bmatrix}R(s)\\\widetilde W(s)\end{bmatrix}\right)=r+x$
and $\begin{bmatrix}R(s)\\\widetilde W(s)\end{bmatrix}$ has 
 $\tilde q_1,\dots, \tilde q_{r+x} $ as  invariant orders at $\infty$, 
 $\bd=(d_1,\dots, d_{n-r-x})$ as  column minimal indices,  and
$\bv=(v_1,\dots, v_{m+z-r-x})$ as row minimal indices.  
 Let $\psi_1(s)$ be the denominator of the first invariant rational function of  $\begin{bmatrix}R(s)\\\widetilde W(s)\end{bmatrix}$,
and let $P(s)=\psi_1(s)R(s)$ and $Q(s)=\begin{bmatrix}P(s)\\\psi_1(s)\widetilde W(s)\end{bmatrix}$. 
By Lemma \ref{lem_polrat},
the polynomial matrix $P(s)\in\FF[s]^{m\times n}$  has $\bc=(c_1,  \dots,  c_{n-r})$  as
	column minimal indices, $\bu=(u_1, \dots,  u_{m-r})$  as  row minimal indices, and its 
invariant orders at $\infty$ are $p_1=\tilde p_1-\deg(\psi_1),\dots, p_r=\tilde p_r-\deg(\psi_1)$. Let $\alpha_1(s)\mid \dots \mid \alpha_r(s)$
be its invariant factors.
Analogously, the polynomial matrix $Q(s)\in\FF[s]^{(m+z)\times n}$  has $\bd=(d_1,  \dots,  d_{n-r-x})$  as 
	column minimal indices, $\bv=(v_1, \dots,  v_{m+z-r-x})$  as  row minimal indices and its 
invariant orders at $\infty$ are $q_1=\tilde q_1-\deg(\psi_1), \dots, q_{r+x}=\tilde q_{r+x}-\deg(\psi_1)$.

Let $A$ be defined as in (\ref{eqbigAioicmirmi}). By Theorem \ref{theoprescrioirmicmipol},
(\ref{eqinterpolioi}) and  (\ref{eqAleq})-(\ref{eqrmimajpolhatsAioi}) hold, where $\hat \ba$ and $\hat \bb$ are defined in (\ref{eqaioicmirmi}) and
(\ref{eqbioicmirmi}), respectively.

From (\ref{eqinterpolioi}) and  (\ref{eqvusAioi})  we obtain (\ref{eqinterratioi}) and (\ref{eqvusAioirat}), respectively.
By Lemma \ref{lemmahx}, from 
(\ref{eqcmimajpolhatsAioi}) we get (\ref{eqdc}).

Let $j\in \{1, \dots,  z-x\}$. Then
$$
\sum_{i=1}^{j}\hat b_i=\sum_{i=1}^{j}\hat{\tilde b}_i
+\min\{0,\sum_{i=r-j+1}^{r}\deg(\alpha_{i})-A\}
\leq \sum_{i=1}^{j}\hat{\tilde b}_i.
$$
Moreover, from (\ref{eqrmimajpolhatsAioi}) and (\ref{eqinterratioi}),
$$
\sum_{i=1}^{z-x}\hat b_i=\sum_{i=1}^{m+z-r-x}  v_i-\sum_{i=1}^{m-r}  u_i=\sum_{i=1}^{z-x}\hat {\tilde b}_i.
$$
Therefore $\hat \bb \prec \hat{\tilde \bb}$, and by Remark \ref{gmajcp} (\ref{eqrmimajpolhatsAioi}) implies  (\ref{eqrmimajrathatsAioi}).

Conversely, 
assume that  (\ref{eqdc}),
(\ref{eqinterratioi}), 
 (\ref{eqvusAioirat}) and (\ref{eqrmimajrathatsAioi})
hold.
Let $\frac{\eta_1(s)}{\varphi_1(s)}$,\dots, $\frac{\eta_r(s)}{\varphi_r(s)}$ be the invariant rational functions of $R(s)$.
Define
$$
\tilde A=\sum_{i=1}^{r+x}\tilde q_i-\sum_{i=1}^{r}\tilde p_i+\sum_{i=1}^{n-r-x}d_i-\sum_{i=1}^{n-r}c_i+\sum_{i=1}^{m+z-r-x}v_i-\sum_{i=1}^{m-r}u_i.
$$
We analyze two cases: $x>0$ and $x=0$. In both cases we  define an integer $Z\geq\deg(\varphi_1)$ and  a monic polynomial $\psi_1(s)$ such that $\varphi_1(s)\mid \psi_1(s)$ and $\deg(\psi_1)=Z$, and  take the polynomial matrix 
$P(s)=\psi_1(s)R(s)$. By Lemma \ref{lem_polrat},  $P(s)$   
has $\bc=(c_1,  \dots,  c_{n-r})$  as
	column minimal indices, $\bu=(u_1, \dots,  u_{m-r})$  as  row minimal indices, its 
invariant orders at $\infty$ are $p_1=\tilde p_1-Z,\dots, p_r=\tilde p_r-Z$, 
and its invariant factors are 
$\alpha_1(s)=\psi_1(s)\frac{\eta_1(s)}{\varphi_1(s)},\dots, \alpha_r(s)=\psi_1(s)\frac{\eta_r(s)}{\varphi_r(s)}$.

Define $q_i=\tilde q_1-Z, \dots, q_{r+x}=\tilde q_{r+x}-Z$ and 
let 
 $A$ be defined as in (\ref{eqbigAioicmirmi}). Then (\ref{eqinterratioi}) and  (\ref{eqinterpolioi}) are equivalent and $A=\tilde A-xZ$.

\begin{itemize}
\item If $x>0$,  define 
$$
\sum_{i=1}^{j}\hat{\tilde a}_i=
  \sum_{i=1}^{m+z-r-x}v_i-\sum_{i=1}^{m-r}u_i+\sum_{i=1}^{r+j}\tilde q_{i+x-j}-\sum_{i=1}^{r}\max\{\tilde p_{i},\tilde q_{i+x-j}\}-\tilde A,
  \quad 1\leq j \leq x,
$$
$$h_j=\min\{i: d_{i-j+1}<c_i\}, \quad 1\leq j \leq x,$$
and let $Z$ be an integer satisfying 
\begin{equation}\label{eqZmax}\begin{array}{rl}
Z\geq &\deg(\varphi_1),\\
xZ\geq &\tilde A,\\
(x-j)Z\geq &\sum_{i=1}^{h_j}c_i-\sum_{i=1}^{h_j-j}d_i-\sum_{i=1}^j\hat {\tilde a}_i, \quad 1\leq j \leq x-1.
\end{array}
\end{equation}
Then,  from   (\ref{eqZmax})    we obtain $A\leq 0$.

We aim to prove that  (\ref{eqAleq})-(\ref{eqrmimajpolhatsAioi}) hold, where $\hat \ba$ and $\hat \bb$ are defined in (\ref{eqaioicmirmi}) and
(\ref{eqbioicmirmi}), respectively.
 Once this is proven, by
 Theorem \ref{theoprescrioirmicmipol} and Remark \ref{remacdecr}-\ref{remitac}, there exists a polynomial matrix $W(s)\in \efe[s]^{z\times n}$ such that  $\rank \left(\begin{bmatrix}P(s)\\W(s)\end{bmatrix}\right)=r+x$
and $\begin{bmatrix}P(s)\\W(s)\end{bmatrix}$  has 
 $q_1,\dots, q_{r+x} $ as  invariant orders at $\infty$, 
 $\bd=(d_1,\dots, d_{n-r-x})$ as  column minimal indices and
$\bv=(v_1,\dots, v_{m+z-r-x})$ as row minimal indices. Let
$\widetilde W(s)=\frac{1}{\psi_1(s)}W(s)$. 
By Lemma \ref{lem_polrat},
 the rational matrix $\begin{bmatrix}R(s)\\\widetilde W(s)\end{bmatrix}$  has the prescribed invariants.

Therefore, we only need to prove that (\ref{eqAleq})-(\ref{eqrmimajpolhatsAioi}) hold.

Since $A\leq 0$, condition (\ref{eqAleq}) holds,  (\ref{eqvusAioirat}) and (\ref{eqvusAioi}) are equivalent, and   $\hat \bb= \hat {\tilde \bb}$, hence (\ref{eqrmimajrathatsAioi}) and (\ref{eqrmimajpolhatsAioi}) are equivalent.

From  (\ref{eqinterpolioi})
we have
$
\sum_{i=1}^x\hat a_i=
\sum_{i=1}^{m+z-r-x}v_i-\sum_{i=1}^{m-r}u_i+\sum_{i=1}^{r+x}q_{i}-\sum_{i=1}^{r}p_{i}- A$. Thus, from (\ref{eqbigAioicmirmi}) we obtain
\begin{equation}\label{eqsumhata}
\sum_{i=1}^x\hat a_i=\sum_{i=1}^{n-r}c_i-\sum_{i=1}^{n-r-x}d_i.
\end{equation}

  From (\ref{eqdc}), $\sum_{i=h_x+1}^{n-r}c_i=\sum_{i=h_x-x+1}^{n-r-x}d_i$ and
\begin{equation}\label{eqhsumhata}
\sum_{i=1}^x\hat a_i=\sum_{i=1}^{h_x}c_i-\sum_{i=1}^{h_x-x}d_i.
\end{equation}

Moreover,
$$ 
  \sum_{i=1}^{j}\hat a_i=\sum_{i=1}^{j}\hat{\tilde a}_i+(x-j)Z, \quad 1\leq j \leq x-1.
$$
From (\ref{eqZmax}) we obtain
\begin{equation}\label{eqcdh}
\sum_{i=1}^{h_j}c_i-\sum_{i=1}^{h_j-j}d_i\leq 
\sum_{i=1}^{j}\hat a_i \quad 1\leq j \leq x-1.
\end{equation}

 From (\ref{eqdc}) we obtain (\ref{eqdcgeq}) (see Lemma \ref{lemmahx}). 
From  (\ref{eqdcgeq}) and  (\ref{eqsumhata})-(\ref{eqcdh}) we get (\ref{eqcmimajpolhatsAioi}).
\item If $x=0$, then $A=\tilde A$ and by   (\ref{eqdc}),   (\ref{eqinterratioi}) 
 (see Remark \ref{remlemmahx})
 and  (\ref{eqvusAioirat}) 
we get
$$
A=\tilde A=\sum_{i=1}^{r}\tilde q_i-\sum_{i=1}^{r}\tilde p_i+\sum_{i=1}^{m+z-r}v_i-\sum_{i=1}^{m-r}u_i\geq 0.
$$
Let  $\tau(s)$ be  a polynomial of $\deg(\tau)= A$ and take $\psi_1(s)=\varphi_1(s)\tau(s)$ and $Z=\deg(\psi_1)=\deg(\varphi_1)+ A$.
Notice that $\alpha_1(s)=\psi_1(s)\frac{\eta_1(s)}{\varphi_1(s)}=\eta_1(s)\tau(s)$.
Define 
$$
 \beta_1(s)=\eta_1(s), \quad \beta_i(s)=\alpha_i(s),\; 2\leq i \leq r.
$$
Then $\beta_1(s)\mid \dots \mid \beta_r(s)$ and (\ref{eqinterif}) is satisfied.

We aim to prove that (\ref{eqcmimajpol})-(\ref{eqdegsumpolioi})
   hold, where $\bb$ is defined in 
(\ref{eqdefbioi}).
 Once this is proven, by
 Corollary \ref{theoprescr4ioi}, there exists a polynomial matrix $W(s)\in \efe[s]^{z\times n}$ such that  $\rank \left(\begin{bmatrix}P(s)\\W(s)\end{bmatrix}\right)=r$
and $\begin{bmatrix}P(s)\\W(s)\end{bmatrix}$  has 
 $q_1,\dots, q_{r} $ as  invariant orders at $\infty$, 
 $\bd=(d_1,\dots, d_{n-r})$ as  column minimal indices and
$\bv=(v_1,\dots, v_{m+z-r})$ as row minimal indices. Let
$\widetilde W(s)=\frac{1}{\psi_1(s)}W(s)$. 
By Lemma \ref{lem_polrat},
the rational matrix $\begin{bmatrix}R(s)\\\widetilde W(s)\end{bmatrix}$  has the prescribed invariants.

Therefore, we only need to prove that (\ref{eqcmimajpol})-(\ref{eqdegsumpolioi}) hold.

Conditions (\ref{eqcmimajpol}) and  (\ref{eqdc}) are equivalent (see Remarks \ref{gmajcp} and  \ref{remlemmahx}). 
From (\ref{eqinterif}), (\ref{eqinterpolioi}) and  (\ref{eqdc})  we obtain
$$
\begin{array}{rl}
 &\sum_{i=1}^{r}\deg(\lcm(\alpha_i,\beta_{i}))+
  \sum_{i=1}^{r}\max\{p_i,q_{i}\}\\
=&\sum_{i=1}^{r}\deg(\alpha_i)+
  \sum_{i=1}^{r}p_i=\sum_{i=1}^{r}\deg(\beta_i)+A+\sum_{i=1}^{r}p_i\\=&
  \sum_{i=1}^{r}\deg(\beta_i)+
\sum_{i=1}^{m+z-r}v_i-\sum_{i=1}^{m-r} u_i+
\sum_{i=1}^{r}
q_i;
  \end{array}
$$
i.e., (\ref{eqdegsumpolioi}) holds.
Let $j\in \{1, \dots, z\}$. Then
$\lcm(\alpha_i,\beta_{i+j})=\lcm(\alpha_i,\alpha_{i+j})=\alpha_{i+j}=\beta_{i+j}$, $1\leq i \leq r-j$, hence
$
\sum_{i=1}^{j}b_i=\sum_{i=1}^{j}\hat{\tilde b}_i.
$
Therefore $\bb=\hat{\tilde \bb}$ and 
(\ref{eqrmimajpol}) and 
(\ref{eqrmimajrathatsAioi}) are equivalent.
\end{itemize}
\hfill $\Box$

The following is an example of how to obtain a rational completion when the infinite and singular structures are prescribed.

\begin{example}
  \label{exainfsingx0}
  (See Example \ref{exalgclosed}). 
   Let $\FF=\RR$,
  $R(s)=
  \begin{bmatrix}0&s&1\\s^2+1&0&0
  \end{bmatrix}
  \in \FF[s]^{2\times3}
  $, $x=0$ and  $z=1$.
 We prescribe
 $$
\tilde q_1=-6, \quad \tilde q_2=-1, \quad \bd=(1), \quad \bv=(5).$$
There is no  $W(s)\in \RR[s]^{1\times 3}$ such that   $\begin{bmatrix}R(s)\\W(s)\end{bmatrix}$ has the prescribed invariants
(see Example \ref{exalgclosed}).

However, by Theorem \ref{theoprescrioirmicmirat}, there  exists $\widetilde{W}(s)\in\RR(s)^{1\times 3}$ such that $\begin{bmatrix}R(s)\\\widetilde W(s)\end{bmatrix}$ has the prescribed invariants. In order to obtain $\widetilde W(s)$, let $\psi_1(s)=s$ and $P(s)=\psi_1(s)R(s)=
  \begin{bmatrix}
  0&s^2&s\\s(s^2+1)&0&0
  \end{bmatrix}\in \RR[s]^{2\times3}
  $. The matrix 
  $P(s)$ has $\alpha_1(s)=s$, $\alpha_2(s)=s(s^2+1)$ as invariant factors and 
  $p_1=-3$, $p_2=-2$ as invariant orders at $\infty$.

  Let $ q_1=-7$,  $q_2=-2$, $\beta_1(s)=1, \beta_2(s)=s(s^2+1)$. The matrix
  $\begin{bmatrix}P(s)\\W(s)\end{bmatrix}=\begin{bmatrix}
0&s^2&s\\s(s^2+1)&0&0\\(s^2+1)(s^5+1)&0&0
  \end{bmatrix}$ has
$\beta_1(s), \beta_2(s)$ as invariant factors, 
  $q_1, q_2$ as invariant orders at $\infty$,  $\bd$ as column minimal indices and  $\bv$ as row minimal indices. Then, by Lemma \ref{lem_polrat}, the rational matrix 
  $\begin{bmatrix}R(s)\\\widetilde W(s)\end{bmatrix}=\frac{1}{s}
\begin{bmatrix}P(s)\\W(s)\end{bmatrix}=
  \begin{bmatrix}
0&s&1\\s^2+1&0&0\\
(s^2+1)(s^4+\frac{1}{s})&0&0
  \end{bmatrix}$
  has the prescribed invariants.
\end{example}

\medskip
In the following two theorems Problem \ref{problemrat} is solved when the infinite structure and the column minimal indices are prescribed, for polynomial matrices first, and then for rational matrices.

The proof of Theorem \ref{corprescrioicmipol} is analogous to that of \cite[Corollary 3.4]{AmBaMaRo24_2}. It can be found in \ref{secappendix}.

\begin{theorem} {\rm (Prescription of the infinite structure and the column minimal indices for polynomial matrices)}
  \label{corprescrioicmipol}
  Let $P(s)\in\efe[s]^{m\times n}$ be a polynomial matrix, $\rank (P(s))=r$.
Let 
$p_1, \dots,  p_r$ be the invariant orders at $\infty$ and 
$\bc=(c_1,  \dots,  c_{n-r})$  the
column minimal indices of $P(s)$.

Let $z$ and $x$ be integers such that $0\leq x\leq \min\{z, n-r\}$ and 
let  $q_1\leq\dots \leq q_{r+x}$ be integers,  and 
 $\bd=(d_1, \dots, d_{n-r-x})$ be   a partition. 
There exists $W(s)\in \efe[s]^{z\times n}$ such that  $\rank \left(\begin{bmatrix}P(s)\\W(s)\end{bmatrix}\right)=r+x$
and $\begin{bmatrix}P(s)\\W(s)\end{bmatrix}$  has 
 $q_1,\dots, q_{r+x} $ as  invariant orders at $\infty$ and
 $\bd=(d_1,\dots, d_{n-r-x})$ as  column minimal indices
if and only if 
(\ref{eqinterpolioi}),
\begin{equation}\label{eqcdioicmi}
 \sum_{i=1}^{n-r}c_i-\sum_{i=1}^{n-r-x}d_i\geq \sum_{i=1}^{r}\max\{p_{i}, q_{i+x}\}+\sum_{i=1}^{x}q_i-\sum_{i=1}^{r}p_i,
\end{equation}
\begin{equation}\label{eqcmimajhatap}  \bc \prec'  (\bd , \hat \ba'),\end{equation}
where $\hat \ba' = (\hat a'_1, \dots,  \hat a'_x )$ is defined as
\begin{equation}\label{eqaioicmi} 
\begin{array}{rl}\sum_{i=1}^{j}\hat a'_i=&
  \sum_{i=1}^{n-r}c_i-\sum_{i=1}^{n-r-x}d_i
  +\sum_{i=1}^{r}p_i-\sum_{i=1}^{x-j}q_i\\&-\sum_{i=1}^{r}\max\{p_{i}, q_{i+x-j}\},\quad 1\leq j \leq x.
\end{array}
\end{equation}
\end{theorem}

\begin{theorem} {\rm (Prescription of the infinite structure and the column minimal indices for rational matrices)}
  \label{corprescrioicmirat}
  Let $R(s)\in\efe(s)^{m\times n}$ be a rational matrix, $\rank (R(s))=r$.
Let
$\tilde p_1, \dots,  \tilde p_r$ be the invariant orders at $\infty$ and 
$\bc=(c_1,  \dots,  c_{n-r})$  the
column minimal indices of $R(s)$.

Let $z$ and $x$ be integers such that $0\leq x\leq \min\{z, n-r\}$ and 
let  $\tilde q_1\leq\dots \leq \tilde q_{r+x}$ be integers,  and 
 $\bd=(d_1, \dots, d_{n-r-x})$ be a partition. 
There exists $\widetilde W(s)\in \efe(s)^{z\times n}$ such that  $\rank \left(\begin{bmatrix}R(s)\\\widetilde W(s)\end{bmatrix}\right)=r+x$
and $\begin{bmatrix}R(s)\\\widetilde W(s)\end{bmatrix}$  has 
 $\tilde q_1,\dots, \tilde q_{r+x}$ as  invariant orders at $\infty$ and
 $\bd=(d_1,\dots, d_{n-r-x})$ as  column minimal indices
if and only if (\ref{eqdc}) and (\ref{eqinterratioi}) 
 hold. 
\end{theorem}
{\noindent\bf Proof.}
Assume that  there exists $\widetilde W(s)\in \efe(s)^{z\times n}$ 
such that
$\rank \left(\begin{bmatrix}R(s)\\\widetilde W(s)\end{bmatrix}\right)=r+x$
and $\begin{bmatrix}R(s)\\\widetilde W(s)\end{bmatrix}$  has 
 $\tilde q_1,\dots, \tilde q_{r+x} $ as  invariant orders at $\infty$ and
 $\bd=(d_1,\dots, d_{n-r-x})$ as  column minimal indices. By Theorem \ref{theoprescrioirmicmirat},
 (\ref{eqdc}) and 
 (\ref{eqinterratioi}) 
hold.

Conversely, assume that (\ref{eqdc}) and 
(\ref{eqinterratioi})  hold.
Let $\bu=(u_1, \dots, u_{m-r})$  be the row minimal indices of $R(s)$ and let
$$
X=\sum_{i=1}^{  r}\max\{{\tilde p}_{i}, {\tilde q}_{i+x}\}-\sum_{i=1}^{r}{\tilde q}_{i+x}.
$$
 Clearly,  $X\geq0$. 
Define
$$
 \sum_{i=1}^{ j} \hat {\tilde b}'_i=X                  
+\sum_{i=1}^{  r-j}\tilde q_{i+x+j}-
\sum_{i=1}^{  r-j}\max\{\tilde p_i, \tilde q_{i+x+j}\},\quad 1\leq j \leq z-x.
$$
Let
$\hat{\tilde\bb}' = (\hat{\tilde b}'_1, \dots, \hat{\tilde b}'_{z-x} )$ and $\bv=\bu \cup \hat  {\tilde\bb}'$. Thus, 
$\sum_{i=1}^{m+z-r-x}v_i-\sum_{i=1}^{m-r}u_i=\sum_{i=1}^{z-x}\hat {\tilde b}'_i$ and from (\ref{eqinterratioi}) we obtain
\begin{equation*}
  \sum_{i=1}^{m+z-r-x}v_i-\sum_{i=1}^{m-r}u_i=X.
\end{equation*}
 Then (\ref{eqvusAioirat}) holds
and $\hat{\tilde \bb}'=\hat{\tilde \bb}$, where  $\hat{\tilde \bb}$ is defined as in (\ref{eqbioicmirmirat}).  Hence  $\bv=\bu \cup \hat {\tilde\bb}$. By Lemma \ref{lemmacup}, (\ref{eqrmimajrathatsAioi}) holds.
By Theorem \ref{theoprescrioirmicmirat}, the sufficiency of  conditions 
(\ref{eqdc}) and (\ref{eqinterratioi})
 follows. 
\hfill $\Box$

\medskip
In the following two theorems Problem \ref{problemrat} is solved when the infinite structure and the row minimal indices are prescribed, first for polynomial matrices and then for rational matrices.
In the polynomial case the field is required to be algebraically closed.

The proof of Theorem \ref{corprescrioirmipol}
is analogous to that of 
\cite[Corollary 3.6]{AmBaMaRo24_2}. It can be found in 
\ref{secappendix}.
\begin{theorem} {\rm (Prescription of the infinite structure and the row minimal indices for polynomial matrices)}
  \label{corprescrioirmipol}
  Let $\efe$ be an algebraically closed field.
  Let $P(s)\in\efe[s]^{m\times n}$ be a polynomial matrix, $\rank (P(s))=r$.
Let $\alpha_1(s)\mid\cdots\mid\alpha_r(s)$ be the invariant factors,
$p_1, \dots,  p_r$ the invariant orders at $\infty$,
$\bc=(c_1,  \dots,  c_{n-r})$  the
column minimal indices, and $ \bu=(u_1, \dots,  u_{m-r})$  the row minimal
indices of $P(s)$.

Let $z$ and $x$ be integers such that $0\leq x\leq \min\{z, n-r\}$ and let $q_1\leq\dots \leq q_{r+x}$ be integers and $\bv=(v_1, \dots, v_{m+z-r-x})$ be  a  partition.
Let
\begin{equation}\label{eqbigAioirmi}
A'=\sum_{i=1}^{r+x}q_i-\sum_{i=1}^{r}p_i-\sum_{i=1}^{x}c_i+\sum_{i=1}^{m+z-r-x}v_i-\sum_{i=1}^{m-r}u_i.
\end{equation}
There exists $W(s)\in \efe[s]^{z\times n}$ such that  $\rank \left(\begin{bmatrix}P(s)\\W(s)\end{bmatrix}\right)=r+x$
and $\begin{bmatrix}P(s)\\W(s)\end{bmatrix}$  has 
 $q_1,\dots, q_{r+x} $ as  invariant orders at $\infty$ and
$\bv=(v_1,\dots, v_{m+z-r-x})$ as row minimal indices
if and only if 
 (\ref{eqinterpolioi}), 
\begin{equation}\label{eqApleq}
A'\leq \sum_{i=r+x-z+1}^{r}\deg(\alpha_i),
\end{equation}
\begin{equation}\label{eqvusApioi}
 \sum_{i=1}^{m+z-r-x}v_i-\sum_{i=1}^
 {m-r}u_i\geq \max\{0, A'\}+\sum_{i=1}^{r}\max\{p_{i}, q_{i+x}\}-\sum_{i=1}^{r}q_{i+x},
\end{equation}
\begin{equation}\label{eqcmimajpolhatsAioirmi}  (c_1, \dots,c_x) \prec \hat \ba'\end{equation}
\begin{equation}\label{eqrmimajpolhatsAioicmi} \bv \prec'  (\bu , \hat \bb'),\end{equation}
where $\hat \ba' = (\hat a'_1, \dots, \hat a'_x )$ and $\hat \bb' = (\hat b'_1, \dots,  \hat b'_{z-x})$ 
are defined as
\begin{equation}\label{eqaioirmi}\begin{array}{rl}
\sum_{i=1}^{j}\hat a'_i=&
  \sum_{i=1}^{m+z-r-x}v_i-\sum_{i=1}^{m-r}u_i+\sum_{i=1}^{r+j}q_{i+x-j}\\&-\sum_{i=1}^{r}\max\{p_{i}, q_{i+x-j}\}-A',
\quad 1\leq j \leq x,
\end{array}
\end{equation}
\begin{equation}\label{eqbioirmi}
\begin{array}{rl}
 \sum_{i=1}^{j} \hat b'_i=&
			\sum_{i=1}^{m+z-r-x}  v_i-\sum_{i=1}^{m-r}  u_i+\min\{0,\sum_{i=r-j+1}^r\deg(\alpha_{i})-A' \}                       
\\&+\sum_{i=1}^{  r-j}q_{i+x+j}-
\sum_{i=1}^{  r-j}\max\{p_i,  q_{i+x+j}\},
  \quad  1\leq j \leq z-x.
		\end{array}
\end{equation}
\end{theorem}

\begin{rem}\label{remac2}
As in Theorem \ref{theoprescrioirmicmipol}, the necessity of the conditions, and the sufficiency when $A'\leq 0$ hold for arbitrary fields.
\end{rem}

\begin{theorem} {\rm (Prescription of the infinite structure and the row minimal indices for  rational matrices)}
  \label{corprescrioirmirat}
  Let $R(s)\in\efe(s)^{m\times n}$ be a rational matrix, $\rank (R(s))=r$. 
 Let  
$\tilde p_1, \dots, \tilde p_r$ be its invariant orders at $\infty$, 
	and $\bu=(u_1, \dots,  u_{m-r})$  its row minimal
	indices.
	
	Let $z,x$ be integers such that $0\leq x\leq \min\{z, n-r\}$ and
        let $\tilde q_{1}\leq \dots \leq \tilde q_{r+x}$ be integers, and $\bv=(v_1, \dots,  v_{m+z-r-x})$  be a partition.
	There exists a rational matrix  $\widetilde W(s)\in \efe(s)^{z\times n}$ such that $\rank \left(\begin{bmatrix}R(s)\\\widetilde W(s)\end{bmatrix}\right)=r+x$ and
        $\begin{bmatrix}R(s)\\\widetilde W(s)\end{bmatrix}$ has 
 $\tilde q_1, \dots,\tilde q_{r+x}$ as invariant orders at $\infty$ and $v_1, \dots,  v_{m+z-r-x}$ as row minimal indices  
        if and only if 
 (\ref{eqinterratioi}),  (\ref{eqvusAioirat}) and (\ref{eqrmimajrathatsAioi}) hold, where 
        $\hat{\tilde \bb}=(\hat{\tilde b}_1, \dots, \hat{\tilde b}_{z-x})$ is defined in (\ref{eqbioicmirmirat}).

\end{theorem}
{\noindent\bf Proof.}
Assume that there exists $\widetilde W(s)\in \efe(s)^{z\times n}$ such that  $\rank \left(\begin{bmatrix}R(s)\\\widetilde W(s)\end{bmatrix}\right)=r+x$
and $\begin{bmatrix}R(s)\\\widetilde W(s)\end{bmatrix}$ has 
 $\tilde q_1,\dots, \tilde q_{r+x} $ as  invariant orders at $\infty$ and
$\bv=(v_1,\dots, v_{m+z-r-x})$ as row minimal indices.  
By 
Theorem \ref{theoprescrioirmicmirat}, 
 (\ref{eqinterratioi}),  (\ref{eqvusAioirat}) and (\ref{eqrmimajrathatsAioi}) hold.

Conversely, assume that 
 (\ref{eqinterratioi}),  (\ref{eqvusAioirat}) and (\ref{eqrmimajrathatsAioi}) hold.
 Let
$\bc=(c_1, \dots, c_{n-r})$ be the column minimal indices of $R(s)$. Define
$\bd=(d_1, \dots, d_{n-r-x})$ as $d_i=c_{i+x}$ for $1\leq i\leq n-r-x$. Then  (\ref{eqdc}) holds and the result follows from Theorem  \ref{theoprescrioirmicmirat}.
\hfill $\Box$

\subsection{Prescription of the finite and singular structures}\label{subsecfinsing}

In this subsection, where the finite structure is prescribed, the results for polynomial matrices are obtained as particular cases of the respective rational theorems.  As in Corollary \ref{theoprescr4ioi}, prescribing $\psi_1(s)=1$, the completed matrices  are polynomial.

We start by prescribing the finite structure and the column and row minimal indices.

\begin{theorem} {\rm (Prescription of the finite and singular structures)}
  \label{theoprescrifrmicmirat}
  Let $R(s)\in\efe(s)^{m\times n}$ be a rational matrix, $\rank (R(s))=r$. 
 Let $\frac{\eta_1(s)}{\varphi_1(s)},\dots, \frac{\eta_r(s)}{\varphi_r(s)}$ be its invariant rational functions,  
	$\bc=(c_1,  \dots,  c_{n-r})$  its
	column minimal indices, and $\bu=(u_1, \dots,  u_{m-r})$  its row minimal
	indices.

Let $z,x$ be integers such that $0\leq x\leq \min\{z, n-r\}$ and let $\epsilon_1(s)\mid \dots \mid\epsilon_{r+x}(s)$ and 
        $\psi_{r+x}(s)\mid \dots \mid\psi_{1}(s)$ be monic polynomials such that $\frac{\epsilon_i(s)}{\psi_i(s)}$ are irreducible rational functions, $1\leq i \leq r+x$, and 
        $\bd=(d_1,  \dots,  d_{n-r-x})$  and $\bv=(v_1, \dots,  v_{m+z-r-x})$  be two partitions.
	There exists a rational matrix  $\widetilde W(s)\in \efe(s)^{z\times n}$ such that $\rank \left(\begin{bmatrix}R(s)\\\widetilde W(s)\end{bmatrix}\right)=r+x$ and
        $\begin{bmatrix}R(s)\\\widetilde W(s)\end{bmatrix}$ has 
$\frac{\epsilon_1(s)}{\psi_1(s)},\dots, \frac{\epsilon_{r+x}(s)}{\psi_{r+x}(s)}$ as invariant rational functions,  
  $d_1,  \dots,  d_{n-r-x}$  as column minimal indices and $v_1, \dots,  v_{m+z-r-x}$ as row minimal indices  
        if and only if  (\ref{eqdc}),
 (\ref{eqinterratnum}), (\ref{eqinterratden}),
 \begin{equation}\label{eqvusifcmirmigrat}
 \sum_{i=1}^{m+z-r-x}v_i-\sum_{i=1}^{m-r}u_i\geq \sum_{i=1}^{r}\Delta\left(\frac{\eta_{i}}{\varphi_{i}},\frac{\epsilon_{i+x}}{\psi_{i+x}}\right)-\sum_{i=1}^{r}\Delta\left(\frac{\epsilon_{i+x}}{\psi_{i+x}}\right),
\end{equation}
\begin{equation}\label{eqrmimajratif} \bv \prec'  (\bu , \hat{\tilde \bb}'),\end{equation}
where $\hat{\tilde \bb}'=( \hat {\tilde b}'_1, \dots,  \hat {\tilde b}'_{z-x})$ is defined as
\begin{equation}\label{eqbifcmirmirat}
\begin{array}{rl}
 \sum_{i=1}^{j} \hat {\tilde b}'_i=&
			\sum_{i=1}^{m+z-r-x}  v_i-\sum_{i=1}^{m-r}  u_i                    
+\sum_{i=1}^{  r-j}\Delta\left(\frac{\epsilon_{i+x+j}}{\psi_{i+x+j}}\right)\\&-\sum_{i=1}^{  r-j}\Delta\left(\frac{\eta_{i}}{\varphi_{i}},\frac{\epsilon_{i+x+j}}{\psi_{i+x+j}}\right), \quad  1\leq j \leq z-x.
		\end{array}
\end{equation}
\end{theorem}
{\noindent\bf Proof.} 
Let $\tilde p_1,\dots, \tilde p_{r} $  be the  invariant orders at $\infty$ of $R(s)$.

Assume that there exists $\widetilde W(s)\in \efe(s)^{z\times n}$ such that  $\rank \left(\begin{bmatrix}R(s)\\\widetilde W(s)\end{bmatrix}\right)=r+x$
and $\begin{bmatrix}R(s)\\\widetilde W(s)\end{bmatrix}$ has
$\frac{\epsilon_1(s)}{\psi_1(s)},\dots, \frac{\epsilon_{r+x}(s)}{\psi_{r+x}(s)}$ as invariant rational functions,  
 $\bd=(d_1,\dots, d_{n-r-x})$ as  column minimal indices,  and
$\bv=(v_1,\dots, v_{m+z-r-x})$ as row minimal indices.  
Let $\tilde q_1,\dots, \tilde q_{r+x} $ be the invariant orders at $\infty$ of $\begin{bmatrix}R(s)\\\widetilde W(s)\end{bmatrix}$.  

By Theorem \ref{prescr4rat}, (\ref{eqinterratnum})-(\ref{eqdegsumrat}) hold, where $\tilde \ba$ and $\tilde \bb$ are defined in (\ref{eqdeftildea}) and (\ref{eqdeftildeb}), respectively.
From (\ref{eqcmimajrat}) and Lemma \ref{lemmahx} we obtain (\ref{eqdc}), and from (\ref{eqdegsumrat}) we obtain (\ref{eqvusifcmirmigrat}).

 We have
$$
\sum_{i=1}^{j}\tilde b_i\leq \sum_{i=1}^{j}\hat{\tilde b}'_i, \quad 1\leq j \leq z-x.
$$
Moreover, from (\ref{eqinterratnum})-(\ref{eqinterratioi}) we obtain 
$$
\sum_{i=1}^{z-x}\tilde b_j=\sum_{i=1}^{z-x}\hat {\tilde b}'_j=\sum_{i=1}^{m+z-r-x}  v_i-\sum_{i=1}^{m-r}  u_i.
$$
Therefore $\tilde \bb \prec \hat{\tilde \bb}'$, and by Remark \ref{gmajcp}  (\ref{eqrmimajrat}) implies (\ref{eqrmimajratif}).

Conversely, assume that (\ref{eqdc}),
(\ref{eqinterratnum}), (\ref{eqinterratden}), (\ref{eqvusifcmirmigrat}) and (\ref{eqrmimajratif})
hold.
We will define integers $\tilde q_1\leq \dots \leq  \tilde q_{r+x}$ such that (\ref{eqinterratioi})-(\ref{eqdegsumrat}) hold, where $\tilde \ba$ and $\tilde \bb$ are defined in
(\ref{eqdeftildea}) and (\ref{eqdeftildeb}), respectively.
Then, by
 Theorem \ref{prescr4rat}, there exists a rational matrix $\widetilde W(s)\in \efe(s)^{z\times n}$ such that  $\rank \left(\begin{bmatrix}R(s)\\\widetilde W(s)\end{bmatrix}\right)=r+x$
and $\begin{bmatrix}R(s)\\\widetilde W(s)\end{bmatrix}$  has 
$\frac{\epsilon_1(s)}{\psi_1(s)},\dots, \frac{\epsilon_{r+x}(s)}{\psi_{r+x}(s)}$ as invariant rational functions,  
  $d_1,  \dots,  d_{n-r-x}$  as column minimal indices, $v_1, \dots,  v_{m+z-r-x}$ as row minimal indices and
 $\tilde q_1,\dots, \tilde q_{r+x} $ as  invariant orders at $\infty$.
 
Therefore, we only need to  define integers $\tilde q_1\leq \dots \leq  \tilde q_{r+x}$ such that (\ref{eqinterratioi})-(\ref{eqdegsumrat}) are satisfied.
Let
$$\tilde B=\sum_{i=1}^{r+x}\Delta\left(\frac{\epsilon_i}{\psi_i}\right)-\sum_{i=1}^{r}\Delta\left(\frac{\eta_i}{\varphi_i}\right)+\sum_{i=1}^{n-r-x}d_i-\sum_{i=1}^{n-r}c_i+\sum_{i=1}^{m+z-r-x}v_i-\sum_{i=1}^{m-r}u_i.$$
\begin{itemize} 
\item If $x>0$, define
$$\begin{array}{rl}
\sum_{i=1}^{j}\hat{\tilde a}'_i=&
  \sum_{i=1}^{m+z-r-x}v_i-\sum_{i=1}^{m-r}u_i+\sum_{i=1}^{r+j}\Delta (\frac{\epsilon_{i+x-j}}{\psi_{i+x-j}})\\&-\sum_{i=1}^{r}
\Delta\left(\frac{\eta_i}{\varphi_i}, \frac{\epsilon_{i+x-j}}{\psi_{i+x-j}}\right),  \quad 1\leq j \leq x,
\end{array}
$$
$$h_0=0, \quad \quad h_j=\min\{i: d_{i-j+1}<c_i\}, \quad 1\leq i \leq x,$$
$$
T_j=\sum_{i=1}^{h_j}c_i-\sum_{i=1}^{h_j-j}d_i-\sum_{i=1}^j\hat {\tilde a}'_i-j\tilde p_1, \quad 0\leq j\leq x,
$$
$$
Z_2=\max\{T_j\; : \; 0\leq j \leq x\}, \quad Z_1=Z_2+\tilde B+x\tilde p_1.
$$
Observe that $Z_2\geq T_0= 0$.
Moreover, from (\ref{eqdc}), (\ref{eqinterratnum}) and (\ref{eqinterratden})  we obtain
\begin{equation}\label{eqsumap}
\sum_{i=1}^{x}\hat{\tilde a}'_i=\tilde B+\sum_{i=1}^{n-r}c_i-\sum_{i=1}^{n-r-x}d_i=
\tilde B+\sum_{i=1}^{h_x}c_i-\sum_{i=1}^{h_x-x}d_i.\end{equation}
Therefore, $T_x=-\tilde B-x\tilde{p}_1$
and  $Z_1=Z_2-T_x\geq 0$.
Define
$$
\begin{array}{rl}
\tilde q_1=&\tilde p_1-Z_1, \\
\tilde q_i=&\tilde p_1, \quad 2\leq i \leq x,\\
\tilde q_{i+x}=&\tilde p_i, \quad  1\leq i \leq r-1,\\
\tilde q_{r+x}=&\tilde p_r+Z_2.
\end{array}
$$
We have
$\tilde q_1\leq \dots \leq  \tilde q_{r+x}$, $\sum_{i=1}^r\tilde p_i-\sum_{i=1}^{r+x}\tilde q_i=Z_1-Z_2-x\tilde{p}_1=\tilde B$, and
$
\tilde q_i\leq \tilde p_i \leq \tilde q_{i+x}$,  for $1\leq i \leq r$;  hence  (\ref{eqinterratioi}) holds.
Moreover,
$$\sum_{i=1}^{r-j}\max\{\tilde p_i, \tilde q_{i+x+j}\}=\sum_{i=1}^{r-j}\tilde q_{i+x+j}, \quad 0\leq j\leq z-x.$$
 Therefore, $\tilde \bb=\hat{\tilde \bb}'$ and  (\ref{eqdegsumrat}) and (\ref{eqrmimajrat}) are respectively equivalent to 
 (\ref{eqvusifcmirmigrat}) and   (\ref{eqrmimajratif}).

Let $j\in\{1, \dots x-1\}$. Then
$$\begin{array}{rl}
\sum_{i=1}^j\tilde a_i=&
\sum_{i=1}^j\hat{\tilde a}'_i+\sum_{i=1}^{r+j}\tilde q_{i+x-j}-\sum_{i=1}^{r}\max\{\tilde p_i, \tilde q_{i+x-j}\}\\
=&\sum_{i=1}^j\hat{\tilde a}'_i+j\tilde p_1+\sum_{i=1}^{r}\tilde p_{i}+Z_2-\sum_{i=1}^j\max\{\tilde p_i,\tilde p_1\}\\&-\sum_{i=j+1}^{r}\max\{\tilde p_i, \tilde p_{i-j}\}
=\sum_{i=1}^j\hat{\tilde a}'_i+j\tilde p_1+Z_2.
\end{array}
$$
From the definition of $Z_2$ we get
\begin{equation}\label{hj}
\sum_{i=1}^{h_j}c_i-\sum_{i=1}^{h_j-j}d_i\leq \sum_{i=1}^j\tilde a_i, \quad 1\leq j \leq x-1.
\end{equation}
Moreover,
$$\begin{array}{rl}
\sum_{i=1}^x\tilde a_i=&
\sum_{i=1}^x\hat{\tilde a}'_i+\sum_{i=1}^{r+x}\tilde q_{i}-\sum_{i=1}^{r}\max\{\tilde p_i, \tilde q_{i}\}\\=&\sum_{i=1}^x\hat{\tilde a}'_i+Z_2-Z_1+x\tilde p_1=
\sum_{i=1}^x\hat{\tilde a}'_i-\tilde B, 
\end{array}$$
and from  (\ref{eqsumap})  we obtain
\begin{equation}\label{eqsuma}
\sum_{i=1}^{x}\tilde a_i=\sum_{i=1}^{n-r}c_i-\sum_{i=1}^{n-r-x}d_i=
\sum_{i=1}^{h_x}c_i-\sum_{i=1}^{h_x-x}d_i.\end{equation}
From (\ref{eqdc}) and Lemma \ref{lemmahx} we get (\ref{eqdcgeq}).
Finally, from  (\ref{eqdcgeq})  and (\ref{eqsuma}) we obtain (\ref{eqcmimajrat}).

\item
If $x=0$, then by (\ref{eqdc}), $\bc=\bd$ (see Remark \ref{remlemmahx}), hence (\ref{eqcmimajrat}) holds and
$$\tilde B=\sum_{i=1}^{r}
\Delta\left(\frac{\epsilon_i}{\psi_i}\right)-\sum_{i=1}^{r}\Delta\left(\frac{\eta_i}{\varphi_i}\right)+\sum_{i=1}^{m+z-r}v_i-\sum_{i=1}^{m-r}u_i.$$
From (\ref{eqinterratnum}), (\ref{eqinterratden}) and
(\ref{eqvusifcmirmigrat}) we obtain $\tilde B\geq 0$.
Define
$$
\begin{array}{rl}
\tilde q_1=&\tilde p_1-\tilde B, \\
\tilde q_i=&\tilde p_i, \quad 2\leq i \leq r.\\
\end{array}
$$
We have
$\tilde q_1\leq \dots \leq  \tilde q_{r}$, $\sum_{i=1}^r\tilde p_i-\sum_{i=1}^{r}\tilde q_i=\tilde B$, and
(\ref{eqinterratioi}) holds.

From (\ref{eqinterratnum})-(\ref{eqinterratioi}) we obtain
$$
\begin{array}{l}
\sum_{i=1}^{m+z-r-x}v_i-\sum_{i=1}^{m-r}u_i
  +\sum_{i=1}^{r}\Delta\left(\frac{\epsilon_i}{\psi_i}\right) +\sum_{i=1}^{r}\tilde q_{i}\\=
  \sum_{i=1}^{m+z-r-x}v_i-\sum_{i=1}^{m-r}u_i
  +\sum_{i=1}^{r}\Delta\left(\frac{\epsilon_i}{\psi_i}\right) +\sum_{i=1}^{r}\tilde p_{i}-\tilde B\\=
 \sum_{i=1}^{r}\Delta\left(\frac{\eta_i}{\varphi_i}\right)+\sum_{i=1}^{r}\tilde p_i=
 \sum_{i=1}^{r}
  \Delta\left(\frac{\eta_i}{\varphi_i}, \frac{\epsilon_{i}}{\psi_{i}}, \tilde p_i,\tilde q_{i}\right);
\end{array}
$$
i.e, (\ref{eqdegsumrat}) holds.

Let $j\in\{1, \dots z-x\}$. Then
$\sum_{i=1}^{r-j}\max\{\tilde p_i, \tilde q_{i+j}\}=\sum_{i=1}^{r-j}\tilde q_{i+j}$; hence $\tilde \bb=\hat{\tilde \bb}'$ and 
(\ref{eqrmimajrat}) is equivalent to (\ref{eqrmimajratif}).

\end{itemize}
\hfill $\Box$

\medskip
The proofs of Theorems \ref{corprescrifcmi} and \ref{corprescrifrmi} follow the scheme of those of Theorems  \ref{corprescrioicmirat} and  \ref{corprescrioirmirat}, respectively, by substituting the roles of the invariant orders at $\infty$ by the invariant rational functions and involving Theorem \ref{theoprescrifrmicmirat} instead of 
Theorem \ref{theoprescrioirmicmirat}.
\begin{theorem} {\rm (Prescription of the finite structure and the  column minimal indices)}
  \label{corprescrifcmi}
  Let $R(s)\in\efe(s)^{m\times n}$ be a rational matrix, $\rank (R(s))=r$. 
 Let $\frac{\eta_1(s)}{\varphi_1(s)},\dots, \frac{\eta_r(s)}{\varphi_r(s)}$ be its invariant rational functions,  
and	$\bc=(c_1,  \dots,  c_{n-r})$  its
	column minimal indices.

Let $z,x$ be integers such that $0\leq x\leq \min\{z, n-r\}$ and let $\epsilon_1(s)\mid \dots \mid\epsilon_{r+x}(s)$ and 
        $\psi_{r+x}(s)\mid \dots \mid\psi_{1}(s)$ be monic polynomials such that $\frac{\epsilon_i(s)}{\psi_i(s)}$ are irreducible rational functions, $1\leq i \leq r+x$, and 
        $\bd=(d_1,  \dots,  d_{n-r-x})$   be a partition.
	There exists a rational matrix  $\widetilde W(s)\in \efe(s)^{z\times n}$ such that $\rank \left(\begin{bmatrix}R(s)\\\widetilde W(s)\end{bmatrix}\right)=r+x$ and
$\begin{bmatrix}R(s)\\\widetilde W(s)\end{bmatrix}$ has 
$\frac{\epsilon_1(s)}{\psi_1(s)},\dots, \frac{\epsilon_{r+x}(s)}{\psi_{r+x}(s)}$ as invariant rational functions and  
  $\bd=(d_1,\dots, d_{n-r-x})$ as column minimal indices 
        if and only if 
 (\ref{eqdc}), (\ref{eqinterratnum}) and (\ref{eqinterratden})
 hold.
\end{theorem}

\begin{theorem} {\rm (Prescription of the finite structure and the row minimal indices)}
  \label{corprescrifrmi}
  Let $R(s)\in\efe(s)^{m\times n}$ be a rational matrix, $\rank (R(s))=r$. 
 Let $\frac{\eta_1(s)}{\varphi_1(s)},\dots, \frac{\eta_r(s)}{\varphi_r(s)}$ be its invariant rational functions  
and $\bu=(u_1, \dots,  u_{m-r})$  its row minimal
	indices.

Let $z,x$ be integers such that $0\leq x\leq \min\{z, n-r\}$ and let $\epsilon_1(s)\mid \dots \mid\epsilon_{r+x}(s)$ and 
        $\psi_{r+x}(s)\mid \dots \mid\psi_{1}(s)$ be monic polynomials such that $\frac{\epsilon_i(s)}{\psi_i(s)}$ are irreducible rational functions, $1\leq i \leq r+x$, and 
$\bv=(v_1, \dots,  v_{m+z-r-x})$  be a partition.
	There exists a rational matrix  $\widetilde W(s)\in \efe(s)^{z\times n}$ such that $\rank \left(\begin{bmatrix}R(s)\\\widetilde W(s)\end{bmatrix}\right)=r+x$ and
$\begin{bmatrix}R(s)\\\widetilde W(s)\end{bmatrix}$ has 
$\frac{\epsilon_1(s)}{\psi_1(s)},\dots, \frac{\epsilon_{r+x}(s)}{\psi_{r+x}(s)}$ as invariant rational functions and $\bv=(v_1,\dots, v_{m+z-r-x})$ as row minimal indices  
        if and only if
 (\ref{eqinterratnum}), (\ref{eqinterratden}), (\ref{eqvusifcmirmigrat}) and (\ref{eqrmimajratif}), where $\hat{\tilde \bb}'$ is defined in (\ref{eqbifcmirmirat}).
\end{theorem}

\subsection{Prescription of the singular structure}\label{subsecsing}

In this subsection we first prescribe the total singular structure, then only the row minimal indices, and finally, only the column minimal indices. In the three cases, the conditions are the same in the rational and polynomial cases.

\begin{theorem} {\rm (Prescription of the singular structure)}
  \label{theoprescrrmicmirat}
  Let $R(s)\in\efe(s)^{m\times n}$ be a rational matrix ($R(s)\in\efe[s]^{m\times n}$ be a polynomial matrix), $\rank (R(s))=r$. 
 Let 
	$\bc=(c_1,  \dots,  c_{n-r})$  be its
	column minimal indices, and $\bu=(u_1, \dots,  u_{m-r})$  its row minimal
	indices.

Let $z,x$ be integers such that $0\leq x\leq \min\{z, n-r\}$ and let 
        $\bd=(d_1,  \dots,  d_{n-r-x})$  and $\bv=(v_1, \dots,  v_{m+z-r-x})$  be two partitions.
	There exists a rational matrix  $\widetilde W(s)\in \efe(s)^{z\times n}$ (a polynomial matrix  $\widetilde W(s)\in \efe[s]^{z\times n}$) such that $\rank \left(\begin{bmatrix}R(s)\\\widetilde W(s)\end{bmatrix}\right)=r+x$ and $\begin{bmatrix}R(s)\\\widetilde W(s)\end{bmatrix}$ has 
  $\bd=(d_1,  \dots,  d_{n-r-x})$  as column minimal indices and $\bv=(v_1, \dots,  v_{m+z-r-x})$ as row minimal indices  
        if and only if (\ref{eqdc}),
\begin{equation}\label{equvgequ}
 \sum_{i=1}^{m+z-r-x}v_i-\sum_{i=1}^{m-r}u_i\geq 0, 
 \end{equation}      
 \begin{equation}\label{eqvmaydif}
 \bv\prec'(\bu, \hat{\tilde \bb}''),
 \end{equation}
 where $\hat{\tilde\bb}''=(\hat{\tilde b}''_1, \dots, \hat{\tilde b}''_{z-x})$ is defined as
 \begin{equation}\label{eqbcmirmir}\sum_{i=1}^{j}\hat{\tilde b}''_i= \sum_{i=1}^{m+z-r-x}v_i-\sum_{i=1}^{m-r}u_i, \quad 1\leq j\leq z-x,\end{equation}
 (i.e., $\hat{\tilde b}''_1=\sum_{i=1}^{m+z-r-x}v_i-\sum_{i=1}^{m-r}u_i$ and $\hat{\tilde b}''_i=0$ for $2\leq i \leq z-x$).
\end{theorem}
{\noindent\bf Proof.}
Let $\frac{\eta_1(s)}{\varphi_1(s)},\dots, \frac{\eta_r(s)}{\varphi_r(s)}$ be the invariant rational functions of $R(s)$. (If $R(s)$ is polynomial, then $\varphi_1(s)=1$).

Assume that there exists a rational matrix  $\widetilde W(s)\in \efe(s)^{z\times n}$   such that $\rank \left(\begin{bmatrix}R(s)\\\widetilde W(s)\end{bmatrix}\right)=r+x$, and
        $\begin{bmatrix}R(s)\\\widetilde W(s)\end{bmatrix}$ has 
$\bd=(d_1,  \dots,  d_{n-r-x})$  as column minimal indices and $\bv=(v_1, \dots,  v_{m+z-r-x})$ as row minimal indices. Let
$\frac{\epsilon_1(s)}{\psi_1(s)},\dots, \frac{\epsilon_{r+x}(s)}{\psi_{r+x}(s)}$ be its invariant rational functions. 
By Theorem \ref{theoprescrifrmicmirat}, (\ref{eqdc}), 
(\ref{eqinterratnum}), 
(\ref{eqinterratden}),
 (\ref{eqvusifcmirmigrat}) and 
(\ref{eqrmimajratif}) hold, where $\hat{\tilde\bb}'$ is defined in
(\ref{eqbifcmirmirat}).
From (\ref{eqvusifcmirmigrat}) we obtain (\ref{equvgequ}). We have
$$
\sum_{i=1}^{j}\hat{\tilde b}'_i\leq\sum_{i=1}^{j}\hat{\tilde b}''_i, \quad 1\leq j \leq z-x, 
$$
and, from (\ref{eqinterratnum}), (\ref{eqinterratden}) we get
$$
\sum_{i=1}^{z-x}\hat{\tilde b}'_i= \sum_{i=1}^{m+z-r-x}v_i-\sum_{i=1}^{m-r}u_i=\sum_{i=1}^{z-x}\hat{\tilde b}''_i.
$$
Therefore $\hat{\tilde \bb}'\prec\hat{\tilde \bb}''$, and by Remark \ref{gmajcp} from
(\ref{eqrmimajratif}) we obtain (\ref{eqvmaydif}).

Conversely, assume that 
(\ref{equvgequ}) and 
 (\ref{eqvmaydif}) hold.
Define 
$$
\begin{array}{rl}
  \epsilon_{i}(s)=&1, \quad 1\leq i \leq x,\\
  \epsilon_{i+x}(s)=&\eta_i(s), \quad 1\leq i \leq r,\\
  \psi_{i}(s)=&\varphi_1(s), \quad 1\leq i \leq x,\\
  \psi_{i+x}(s)=&\varphi_i(s), \quad 1\leq i \leq r.
  \end{array}
$$
Then
$\epsilon_1(s)\mid \dots \mid\epsilon_{r+x}(s)$ and 
$\psi_{r+x}(s)\mid \dots \mid\psi_{1}(s)$ are monic polynomials such that $\frac{\epsilon_i(s)}{\psi_i(s)}$ are irreducible rational functions, $1\leq i \leq r+x$. We have
$$\epsilon_{i}(s)\mid\eta_i(s)\mid\epsilon_{i+x}(s), \quad
\psi_{i+x}(s)\mid \varphi_i(s)\mid \psi_{i}(s), \quad 1\leq i \leq r;$$
hence (\ref{eqinterratnum}) and (\ref{eqinterratden}) hold. 
Let $\hat{\tilde \bb}'=( \hat {\tilde b}'_1, \dots,  \hat {\tilde b}'_{z-x})$ be defined as in (\ref{eqbifcmirmirat}).
 As
$$
\sum_{i=1}^{r-j}\Delta\left(\frac{\eta_{i}}{\varphi_{i}},\frac{\epsilon_{i+x+j}}{\psi_{i+x+j}}\right)=
\sum_{i=1}^{r-j}\Delta\left(\frac{\epsilon_{i+x+j}}{\psi_{i+x+j}}\right), \quad 0\leq j \leq z-x,
$$
from
 (\ref{equvgequ})  we obtain (\ref{eqvusifcmirmigrat}), 
$\hat{\tilde \bb}'=\hat{\tilde \bb}''$, and from 
(\ref{eqvmaydif}) we get
(\ref{eqrmimajratif}).
 By Theorem \ref{theoprescrifrmicmirat} the sufficiency of conditions (\ref{eqdc}),
(\ref{equvgequ}) and 
(\ref{eqvmaydif}) follows.

Observe that if $R(s)$ is polynomial,  then $\psi_{1}(s)=\varphi_1(s)=1$ and $\widetilde W(s)$ is also polynomial.
\hfill $\Box$

\begin{theorem} {\rm (Prescription of the  row minimal indices)}
  \label{corprescrrmi}
  Let $R(s)\in\efe(s)^{m\times n}$ be a rational matrix ($R(s)\in\efe[s]^{m\times n}$ be a polynomial matrix), $\rank (R(s))=r$. 
 Let 
 $\bu=(u_1, \dots,  u_{m-r})$  be its row minimal
	indices.

Let $z,x$ be integers such that $0\leq x\leq \min\{z, n-r\}$ and let 
        $\bv=(v_1, \dots,  v_{m+z-r-x})$  be a partition.
	There exists a rational matrix  $\widetilde W(s)\in \efe(s)^{z\times n}$ (a polynomial matrix  $\widetilde W(s)\in \efe[s]^{z\times n}$) such that $\rank \left(\begin{bmatrix}R(s)\\\widetilde W(s)\end{bmatrix}\right)=r+x$ and
        $\begin{bmatrix}R(s)\\\widetilde W(s)\end{bmatrix}$ has 
   $\bv=(v_1, \dots,  v_{m+z-r-x})$ as row minimal indices  
        if and only if (\ref{equvgequ}) and 
(\ref{eqvmaydif}) hold, where
        $\hat{\tilde\bb}''=(\hat{\tilde b}''_1, \dots, \hat{\tilde b}''_{z-x})$ is defined in (\ref{eqbcmirmir}).
 
\end{theorem}
{\noindent\bf Proof.}
Assume that there exists a rational matrix  $\widetilde W(s)\in \efe(s)^{z\times n}$  such that $\rank \left(\begin{bmatrix}R(s)\\\widetilde W(s)\end{bmatrix}\right)=r+x$ and
        $\begin{bmatrix}R(s)\\\widetilde W(s)\end{bmatrix}$ has 
$\bv=(v_1, \dots,  v_{m+z-r-x})$ as row minimal indices.
By Theorem \ref{theoprescrrmicmirat},
(\ref{equvgequ}) and 
(\ref{eqvmaydif}) hold.

Conversely, assume that 
(\ref{equvgequ}) and 
(\ref{eqvmaydif}) hold. Let
$\bc=(c_1, \dots, c_{n-r})$ be the column minimal indices of $R(s)$. Define
$\bd=(d_1, \dots, d_{n-r-x})$  as $d_i=c_{i+x}$, for $1\leq i \leq n-r-x$. Then
(\ref{eqdc}) holds. By Theorem \ref{theoprescrrmicmirat} the sufficiency of conditions 
(\ref{equvgequ}) and 
(\ref{eqvmaydif}) follows.
\hfill $\Box$

\begin{theorem} {\rm (Prescription of the column minimal indices)}
  \label{corprescrcmi}
  Let $R(s)\in\efe(s)^{m\times n}$ be a rational matrix ($R(s)\in\efe[s]^{m\times n}$ be a polynomial matrix), $\rank (R(s))=r$. 
 Let 
 $\bc=(c_1,  \dots,  c_{n-r})$  be its column minimal indices.

Let $z,x$ be integers such that $0\leq x\leq \min\{z, n-r\}$ and let 
        $\bd=(d_1,  \dots,  d_{n-r-x})$  be a partition.
	There exists a rational matrix $\widetilde W(s)\in \efe(s)^{z\times n}$ (a polynomial matrix  $\widetilde W(s)\in \efe[s]^{z\times n}$) such that $\rank \left(\begin{bmatrix}R(s)\\\widetilde W(s)\end{bmatrix}\right)=r+x$ and
        $\begin{bmatrix}R(s)\\\widetilde W(s)\end{bmatrix}$ has 
   $\bd=(d_1, \dots, d_{n-r-x})$ as column minimal indices if and only if (\ref{eqdc}) holds.
\end{theorem}
{\noindent\bf Proof.}
Assume that there exists a rational matrix  $\widetilde W(s)\in \efe(s)^{z\times n}$   such that $\rank \left(\begin{bmatrix}R(s)\\\widetilde W(s)\end{bmatrix}\right)=r+x$ and
        $\begin{bmatrix}R(s)\\\widetilde W(s)\end{bmatrix}$ has 
$\bd=(d_1, \dots, d_{n-r-x})$ as column minimal indices.
By Theorem \ref{theoprescrrmicmirat},
(\ref{eqdc}) holds.

Conversely, assume that (\ref{eqdc}) holds.
Let $\bu=(u_1, \dots,  u_{m-r})$  be the row minimal indices of $R(s)$.
Define $\bv=(v_1, \dots, v_{m+z-r-x})$ as
$$\begin{array}{rl}
v_i=&u_i, \quad 1\leq i \leq m-r,\\
v_i=&0, \quad m-r+ 1\leq i \leq m-r+z-x.\\
\end{array}
$$
Obviously,  (\ref{equvgequ}) holds.
Let $\hat{\tilde\bb}''=(\hat{\tilde b}''_1, \dots, \hat{\tilde b}''_{z-x})$   be defined as  in (\ref{eqbcmirmir}). Then $\hat{\tilde\bb}''_i=0$ for $1\leq i \leq z-x$  and 
 $\bv=\bu \cup \hat {\tilde \bb}''$.
By Lemma \ref{lemmacup}, (\ref{eqvmaydif}) holds.
By Theorem \ref{theoprescrifrmicmirat}, the sufficiency of  condition (\ref{eqdc})
 follows.
 \hfill $\Box$

\appendix \section{Proofs}
\label{secappendix}

As announced, we prove Theorems \ref{theoprescrioirmicmipol}, \ref{corprescrioicmipol} and \ref{corprescrioirmipol} in this appendix.

\medskip
\noindent
{\bf Proof of Theorem \ref{theoprescrioirmicmipol}.}
Assume that there exists $W(s)\in \efe[s]^{z\times n}$ such that  $\rank \left(\begin{bmatrix}P(s)\\W(s)\end{bmatrix}\right)=r+x$
and $\begin{bmatrix}P(s)\\W(s)\end{bmatrix}$  has 
 $q_1,\dots, q_{r+x} $ as  invariant orders at $\infty$, 
 $\bd=(d_1,\dots, d_{n-r-x})$ as  column minimal indices and
$\bv=(v_1,\dots, v_{m+z-r-x})$ as row minimal indices. Let $\beta_1(s)\mid \dots \mid \beta_{r+x}(s)$ be its invariant factors.

Then, by Theorem \ref{theoexistencerat} we obtain
\begin{equation}\label{eqbetaIST}
A=\sum_{i=1}^{r}\deg(\alpha_i)-\sum_{i=1}^{r+x}\deg(\beta_i),
\end{equation}
and
by Corollary \ref{theoprescr4ioi} conditions (\ref{eqinterif})-(\ref{eqdegsumpolioi}) hold, where $\ba$ and $\bb$ are defined in (\ref{eqdefaioi}) and (\ref{eqdefbioi}), respectively.
Taking into account that
$$
\begin{array}{rl}
&\sum_{i=1}^{  r}\deg(\lcm(  \alpha_{i},  \beta_{i+x}))-\sum_{i=1}^{r}\deg( \beta_{i+x})
\\\geq&
\max \{\sum_{i=1}^{r}\deg( \alpha_{i})-\sum_{i=1}^{r}\deg( \beta_{i+x}), 0\}\geq \max\{A,0\},\\
\end{array}
$$
from (\ref{eqdegsumpolioi}) we obtain (\ref{eqvusAioi}).
From (\ref{eqinterif}) we have $\sum_{i=1}^{r+x}\deg(\beta_i)=\sum_{i=1}^{z}\deg(\beta_i)+\sum_{i=1}^{r+x-z}\deg(\beta_{i+z})\geq \sum_{i=1}^{r+x-z}\deg(\alpha_{i});$
hence we obtain (\ref{eqAleq}). 
We also have
$$\begin{array}{rl}\sum_{i=1}^{x}\hat a_i= &\sum_{i=1}^{m+z-r-x}v_i-\sum_{i=1}^{m-r}u_i+\sum_{i=1}^{r+x}q_i-\sum_{i=1}^{r}p_i-A\\=&\sum_{i=1}^{n-r}c_i-
\sum_{i=1}^{n-r-x}d_i=\sum_{i=1}^{x}a_i,\end{array}$$
$$
\begin{array}{rl}
  \sum_{i=1}^{j}a_i=&\sum_{i=1}^j \hat a_i+A+\sum_{i=1}^{r+j}\deg(\beta_{i+x-j})-\sum_{i=1}^{r}\deg(\lcm( \alpha_i,  \beta_{i+x-j}))
  \\=&
  \sum_{i=1}^j \hat a_i+\sum_{i=1}^{r}\deg(\alpha_i)-\sum_{i=1}^{x-j}\deg(\beta_i)-\sum_{i=1}^{r}\deg(\lcm( \alpha_i,  \beta_{i+x-j}))\\\leq&
  \sum_{i=1}^j \hat a_i, \quad 1\leq j\leq x-1.
\end{array}
$$
 Therefore $\ba\prec \hat \ba$, and  from (\ref{eqcmimajpol})  and Remark \ref{gmajcp} 
 we derive  (\ref{eqcmimajpolhatsAioi}).
Moreover, from (\ref{eqinterif}), (\ref{eqinterpolioi}) and (\ref{eqAleq}) we get
$\sum_{i=1}^{z-x}\hat b_i=\sum_{i=1}^{m+z-r-x}v_i-\sum_{i=1}^{m-r}u_i=
\sum_{i=1}^{z-x} b_i, $
 and since 
 $$
		\begin{array}{rl}
\sum_{i=1}^{j} b_i=&\sum_{i=1}^{j} \hat b_i-\min\{0,\sum_{i=r-j+1}^{r}\deg(\alpha_{i})-A\}  \\&
                        +\sum_{i=1}^{r-j}\deg(\beta_{i+x+j})
			-\sum_{i=1}^{r-j}\deg(\lcm(\alpha_i, \beta_{i+x+j}))\\
   =&\sum_{i=1}^{j} \hat b_i
   -\sum_{i=1}^{r-j}\deg(\lcm(\alpha_i, \beta_{i+x+j}))\\&+
   \max\{\sum_{i=1}^{r-j}\deg(\beta_{i+x+j}), 
   \sum_{i=1}^{r-j}\deg(\alpha_{i})-\sum_{i=1}^{x+\min\{j,r\}}\deg(\beta_{i})\}
   \\\leq &\sum_{i=1}^{j} \hat b_i, \quad 1\leq j \leq z-x,
   \end{array}$$
we have $\bb\prec \hat \bb$. Thus,    from (\ref{eqrmimajpol})  and Remark \ref{gmajcp}  
we derive  (\ref{eqrmimajpolhatsAioi}).

Conversely, 
assume that 
(\ref{eqinterpolioi}) and (\ref{eqAleq})-(\ref{eqrmimajpolhatsAioi}) hold.
If we prove that there exist monic polynomials $\beta_1(s)\mid\cdots\mid\beta_{r+x}(s)$ such that
(\ref{eqinterif}) and 
(\ref{eqcmimajpol})-(\ref{eqdegsumpolioi}) are satisfied, then by Corollary \ref{theoprescr4ioi}, there exists  $W(s)\in \efe[s]^{z\times n}$
 such that  $\rank \left(\begin{bmatrix}P(s)\\W(s)\end{bmatrix}\right) =r+x$ and 
	$\begin{bmatrix}P(s)\\W(s)\end{bmatrix}$ has $\beta_1(s), \dots, \beta_{r+x}(s)$ as invariant factors,
	$q_1, \dots, q_{r+x}$ as invariant orders at $\infty$,
 $d_1, \dots, d_{n-r-x}$ as column minimal indices
	and 
	$v_1, \dots, v_{m+z-r-x}$ as row minimal indices.

Therefore, we only need to prove the existence of   polynomials $\beta_1(s)\mid\cdots\mid\beta_{r+x}(s)$   satisfying 
the conditions mentioned above.

Observe that if $x=0$,  then from (\ref{eqcmimajpolhatsAioi}) we get $\bc =\bd$ and $A=\sum_{i=1}^{r}q_i-\sum_{i=1}^{r}p_i+\sum_{i=1}^{m+z-r}v_i-\sum_{i=1}^{m-r}u_i$. From
(\ref{eqvusAioi}) we conclude that  $A\geq 0$.
If $x=z$,  then from (\ref{eqAleq}) we obtain $A\leq 0$.  We distinguish two cases. 

\begin{itemize}
  \item  
Let $A>0$.  We saw in Remark \ref{remacdecr}-2 that if $g=\min\{k\geq 0: A\leq \sum_{i=r-k+1}^{r}\deg(\alpha_{i}) \}$ and    
 $A>0$, then $1\leq g\leq z-x$. Therefore  (\ref{eqtarte}) holds, i.e.,  
$0<A-\sum_{i=r-g+2}^{r}\deg(\alpha_{i})\leq \deg(\alpha_{r-g+1})
$.

Let $h=\min\{k: A-\sum_{i=r-g+2}^{r}\deg(\alpha_{i})\leq \deg(\alpha_{k-g+1}) \}$. Obviously, $h\leq r$.  Since 
$\deg(\alpha_{j-g+1})=0< A-\sum_{i=r-g+2}^{r}\deg(\alpha_{i})$ for $j\leq g-1$, we have
$h\geq g$ and
\begin{equation}\label{eqdeg}\deg(\alpha_{h-g})< A-\sum_{i=r-g+2}^{r}\deg(\alpha_{i})\leq \deg(\alpha_{h-g+1}).\end{equation}

    Taking $$w=
    \deg(\alpha_{h-g})+\deg(\alpha_{h-g+1})+\sum_{i=r-g+2}^{r}\deg(\alpha_{i})-A,$$
     we get
    $\deg(\alpha_{h-g})\leq w<\deg(\alpha_{h-g+1})$.
    As $\efe$ is an algebraically closed field, there exists a monic polynomial
    $\tau(s)$  such that
$$\alpha_{h-g}(s)\mid \tau(s)\mid \alpha_{h-g+1}(s), \quad \deg(\tau)=w.$$
Define 
$$
\begin{array}{rll}
  \beta_i(s)=&\alpha_{i-x-g}(s),&1\leq i\leq h+x-1,
  \\
  \beta_{h+x}(s)=&\tau(s),\\
  \beta_i(s)=&\alpha_{i-x-g+1}(s),&h+x+1\leq i\leq r+x.
\end{array}
$$
We have
$\beta_1(s)\mid \dots \mid\beta_{r+x}(s)$, and 
$$\sum_{i=1}^{r+x}\deg(\beta_{i})=\sum_{i=1}^{h-g+1}\deg(\alpha_{i})+\sum_{i=r-g+2}^{r}\deg(\alpha_{i})-A
+\sum_{i=h-g+2}^{r-g+1}\deg(\alpha_{i}),$$ 
 i.e., (\ref{eqbetaIST}) holds. 
Moreover,
$\beta_i(s)\mid \alpha_{i-x-g+1}(s)\mid\alpha_{i}(s)\mid\beta_{i+x+g}(s)\mid\beta_{i+z}(s)$ for $1\leq i \leq r$, 
therefore (\ref{eqinterif}) holds.
Since $\beta_{i+x-j}(s)\mid \alpha_{i-g-j+1}(s)\mid \alpha_{i}(s)$, $1\leq i \leq r$, $0\leq j \leq x$,  it follows that
$$\sum_{i=1}^{r}\deg(\lcm(\alpha_{i}, \beta_{i+x-j}))=\sum_{i=1}^{r}\deg(\alpha_{i}), \quad 0\leq j \leq x.$$
As $h\geq g \geq 1$, $\beta_i(s)=\alpha_{i-x-g}(s)=1$ for $1\leq i \leq x$ and 
$$\sum_{i=1}^{r+x}\deg(\beta_i)=\sum_{i=1}^{r+j}\deg(\beta_{i+x-j}), \quad 0\leq j \leq x.$$ 
Then we can write
$$A=\sum_{i=1}
^{r}\deg(\alpha_{i})-\sum_{i=1}^{r+x}\deg(\beta_i)=\sum_{i=1}^{r}\deg(\lcm(\alpha_{i}, \beta_{i+x-j}))-\sum_{i=1}^{r+j}\deg(\beta_{i+x-j}), \quad 0\leq j\leq x.$$
From this expression for $j=0$ and taking into account (\ref{eqvusAioi}), we get
$$ 
\sum_{i=1}^{r}\deg(\lcm(\alpha_{i}, \beta_{i+x}))+
\sum_{i=1}^{r}\max\{p_i, q_{i+x}\}\leq 
\sum_{i=1}^{r}
\deg(\beta_{i+x})+\sum_{i=1}^{m+z-r-x}v_i-\sum_{i=1}^{m-r}u_i+\sum_{i=1}^{r}q_{i+x}.
$$
In particular, if $x=0$  we obtain
$$\sum_{i=1}^{r}\deg(\lcm(\alpha_{i}, \beta_{i}))-\sum_{i=1}^{r}\deg(\beta_{i}))=A=\sum_{i=1}^{m+z-r}v_i-\sum_{i=1}^{m-r}u_i+\sum_{i=1}^{r}q_{i}-\sum_{i=1}^{r}\max\{p_i, q_i\},$$
hence condition (\ref{eqdegsumpolioi}) is satisfied.
Furthermore, if $\ba$ is defined as in (\ref{eqdefaioi}), then
$$
  \sum_{i=1}^{j}a_i=
\sum_{i=1}^{j}\hat a_i, \quad 1\leq j\leq x,
$$
 therefore $\ba=\hat \ba$ and (\ref{eqcmimajpolhatsAioi}) is equivalent to (\ref{eqcmimajpol}).

Let $\bb$ be as in (\ref{eqdefbioi}). Then
$$\begin{array}{rl}
    \sum_{i=1}^{j}b_i=&\sum_{i=1}^{j}\hat b_i-\min\{0, \sum_{i=r-j+1}^r\deg(\alpha_i)-A\}+ \sum_{i=1}^{  r-j}\deg(\beta_{i+x+j})\\&-\sum_{i=1}^{  r-j}\deg(\lcm(  \alpha_{i},  \beta_{i+x+j})), \quad 1\leq j \leq z-x.
\end{array}$$
Let $j\in\{1,\dots, g-1\}$. Then $\sum_{i=r-j+1}^{r}\deg(\alpha_{i})
<A$  and $\beta_{i+x+j}(s)\mid\alpha_{i+j-g+1}(s)\mid\alpha_{i}(s)$,
$1\leq i\leq r-j$; hence
$$\sum_{i=1}^{  r-j}\deg(\lcm(  \alpha_{i},  \beta_{i+x+j}))=\sum_{i=1}^{r-j}\deg(\alpha_i),$$ and, taking into account (\ref{eqbetaIST}),
and that for $1\leq i \leq x+j$, $\beta_{i}(s)=\alpha_{i-x-g}(s)=1$, 
we obtain 
$$
\begin{array}{rl}
\sum_{i=1}^{j}b_i=&\sum_{i=1}^{j}\hat b_i+A-\sum_{i=r-j+1}^r\deg(\alpha_i)+
\sum_{i=1}^{  r+x}\deg(  \beta_{i})
-\sum_{i=1}^{r-j}\deg(\alpha_i)\\=&
\sum_{i=1}^{j}\hat b_i.
\end{array}
$$
Let $j\in\{g,\dots, z-x\}$. Then $A\leq \sum_{i=r-j+1}^{r}\deg(\alpha_{i})$ and 
$$\alpha_{i}(s)\mid\beta_{i+x+g}(s)\mid\beta_{i+x+j}(s),\quad 1\leq i\leq r-j;$$
hence
$$\sum_{i=1}^{  r-j}\deg(\lcm(  \alpha_{i},  \beta_{i+x+j}))=\sum_{i=1}^{  r-j}\deg(\beta_{i+x+j}),$$ and
$$
    \sum_{i=1}^{j}b_i=\sum_{i=1}^{j}\hat b_i+\sum_{i=1}^{  r-j}\deg(\beta_{i+x+j})-
    \sum_{i=1}^{  r-j}\deg(\beta_{i+x+j})=\sum_{i=1}^{j}\hat b_i.
$$
Therefore $\bb=\hat \bb$ and (\ref{eqrmimajpolhatsAioi}) is equivalent to (\ref{eqrmimajpol}).

\item  
  If  $A\leq 0$, 
  let $\tau(s)$ be a monic polynomial such that $\deg(\tau)=-A$.
  Define
$$
\begin{array}{rll}
  \beta_i(s)=&\alpha_{i-x}(s),&1\leq i\leq r+x-1,
  \\
  \beta_{r+x}(s)=&\alpha_{r}(s)\tau(s).\\
\end{array}
$$
Then $\beta_1(s)\mid \dots \mid \beta_{r+x}(s)$ and  (\ref{eqbetaIST}) holds. 

 If $x=0$, recall that $A\geq 0$. Then, $A=0$, $\tau(s)=1$, 
$\beta_i(s)=\alpha_i(s)$, $1\leq i \leq r$, and (\ref{eqinterif}) holds.

If $x\geq 1$ we have
$$\beta_i(s)=\alpha_{i-x}(s)\mid \alpha_i(s), \quad 1\leq i \leq r,$$
$$\alpha_{i}(s)\mid\alpha_{i+z-x}(s)=\beta_{i+z}(s), \quad 1\leq i \leq r+x-z-1,$$
$$\alpha_{r+x-z}(s)\mid\alpha_{r}(s)\mid \alpha_{r}(s)\tau(s)=\beta_{r+x}(s),$$
$$\alpha_{i}(s)\mid\beta_{i+z}(s)=0, \quad r+x-z< i \leq r,$$
 therefore (\ref{eqinterif}) also holds.

 By definition, 
 $\alpha_{i}(s)\mid \beta_{i+x}(s)$, $1\leq i \leq r$, then
$\sum_{i=1}^{r}\deg(\lcm(\alpha_i,\beta_{i+x}))=\sum_{i=1}^{r}\deg(\beta_{i+x}).$

 If $x=0$, then $0=A=\sum_{i=1}^{r}q_i-\sum_{i=1}^{r}p_i+\sum_{i=1}^{m+z-r}v_i-\sum_{i=1}^{m-r}u_i$ and from (\ref{eqinterpolioi}) we obtain
$\sum_{i=1}^{m+z-r}v_i-\sum_{i=1}^{m-r}u_i=\sum_{i=1}^{r}p_i-\sum_{i=1}^{r}q_i=\sum_{i=1}^{r}\max\{p_i,q_{i}\}-\sum_{i=1}^{r}q_{i}$.
Thus, from (\ref{eqvusAioi}) we obtain (\ref{eqdegsumpolioi}).

Let
$\ba$ and  $\bb$ be  defined as in (\ref{eqdefaioi}) and (\ref{eqdefbioi}), respectively.

Let
$j\in\{1,\dots, x \}$. Then $\beta_{i+x-j}(s)=\alpha_{i-j}(s)\mid\alpha_{i}(s)$, $1\leq i \leq r$; hence 
$\sum_{i=1}^{  r}\deg(\alpha_{i})
=\sum_{i=1}^{  r}\deg(\lcm(\alpha_{i}, \beta_{i+x-j}))$
 and, taking into account (\ref{eqbetaIST}) and that, for $1\leq i \leq x$, $\beta_i(s)=\alpha_{i-x}(s)=1$, we obtain 
$$\begin{array}{rl}
  \sum_{i=1}^{j}a_i=&\sum_{i=1}^{j}\hat a_i+A+\sum_{i=1}^{r+j}\deg(\beta_{i+x-j})-
\sum_{i=1}^{  r}\deg(\lcm(\alpha_{i}, \beta_{i+x-j}))\\=&
\sum_{i=1}^{j}\hat a_i+A+\sum_{i=1}^{r+x}\deg(\beta_{i})-
\sum_{i=1}^{  r}\deg(\alpha_{i})=
\sum_{i=1}^{j}\hat a_i.
\end{array}
$$
Therefore $\ba=\hat \ba$ and (\ref{eqcmimajpolhatsAioi}) is equivalent to (\ref{eqcmimajpol}).

Let $j\in\{1,\dots, z-x\}$.
Then
$\alpha_{i}(s)\mid\alpha_{i+j}(s)\mid\beta_{i+x+j}(s)$, $1\leq i\leq r-j$; hence
$$\sum_{i=1}^{  r-j}\deg(\lcm(  \alpha_{i},  \beta_{i+x+j}))=\sum_{i=1}^{  r-j}\deg(\beta_{i+x+j})
$$ and
   $$	
\sum_{i=1}^{j} b_i=\sum_{i=1}^{j} \hat b_i
                        +\sum_{i=1}^{r-j}\deg(\beta_{i+x+j})
			-\sum_{i=1}^{r-j}\deg(\lcm(\alpha_i, \beta_{i+x+j}))\\
   =\sum_{i=1}^{j} \hat b_i
   $$
Therefore $\bb=\hat \bb$ and (\ref{eqrmimajpolhatsAioi}) is equivalent to (\ref{eqrmimajpol}).
\end{itemize}
\hfill $\Box$

\noindent
{\bf Proof of Theorem  \ref{corprescrioicmipol}.}
Let $\bu=(u_1, \dots, u_{m-r})$  be the row minimal indices of $P(s)$.

Assume that  there exists $W(s)\in \efe[s]^{z\times n}$ 
such that $\rank\left(\begin{bmatrix}P(s)\\W(s)\end{bmatrix}\right)=r+x$ 
and $\begin{bmatrix}P(s)\\W(s)\end{bmatrix}$ has
 $q_1\leq\dots \leq q_{r+x} $ as  invariant orders at $\infty$ and 
 $\bd=(d_1,\dots, d_{n-r-x})$ as column minimal indices. Let
$\bv=(v_1,\dots, v_{m+z-r-x})$  be its row minimal indices. By Theorem \ref{theoprescrioirmicmipol} and Remark \ref{remacdecr}-\ref{remitac}, 
(\ref{eqinterpolioi}), (\ref{eqvusAioi}) and (\ref{eqcmimajpolhatsAioi}) hold, where
$A$ and  $\hat \ba$ are defined in (\ref{eqbigAioicmirmi}) and  (\ref{eqaioicmirmi}), respectively.
 As 
$$
\sum_{i=1}^{m+z-r-x}v_i-\sum_{i=1}^{m-r}u_i+ \sum_{i=1}^{r+x}q_i=A+\sum_{i=1}^{n-r}c_i-\sum_{i=1}^{n-r-x}d_i+\sum_{i=1}^{r}p_i,
$$
we have  
$$
\begin{array}{rl}\sum_{i=1}^{j}\hat a'_i=&
\sum_{i=1}^{m+z-r-x}v_i-\sum_{i=1}^{m-x}u_i-A+\sum_{i=1}^{r+j}q_{i+x-j}\\&
-\sum_{i=1}^{r}\max\{p_{i}, q_{i+x-j}\}
=\sum_{i=1}^{j}\hat a_i,\quad 1\leq j \leq x,
\end{array}
$$
i.e.,
$\hat \ba=\hat\ba'$, and from  (\ref{eqvusAioi}) and  (\ref{eqcmimajpolhatsAioi}) we obtain (\ref{eqcdioicmi}) and  (\ref{eqcmimajhatap}), respectively.

Conversely, assume that 
(\ref{eqinterpolioi}),
(\ref{eqcdioicmi}) and
 (\ref{eqcmimajhatap}) hold.  
Define
$\hat\bb' = (\hat b'_1, \dots, \hat b'_{z-x} )$ as
$$\begin{array}{rl}\sum_{i=1}^{j}\hat b'_i=&
\sum_{i=1}^{n-r}c_i-\sum_{i=1}^{n-r-x}d_i+\sum_{i=1}^{r}p_i-\sum_{i=1}^{x+\min\{j, r\}}q_i\\&
-\sum_{i=1}^{  r-j}\max\{p_{i}, q_{i+x+j}\}, \quad
                        1\leq j \leq  z-x.
		\end{array}
$$
Take
$\bv=\bu \cup \hat \bb'$. 
Then
$\sum_{i=1}^{m+z-r-x}v_i-\sum_{i=1}^{m-r}u_i=\sum_{i=1}^{z-x}\hat b'_i$.
Moreover, from (\ref{eqinterpolioi}) we obtain
\begin{equation}\label{sumvu}
\sum_{i=1}^{m+z-r-x}v_i-\sum_{i=1}^{m-r}u_i=\sum_{i=1}^{n-r}c_i-\sum_{i=1}^{n-r-x}d_i+\sum_{i=1}^{r}p_i-\sum_{i=1}^{  r+x}q_i.
\end{equation}
Let $\alpha_1(s)\mid\cdots\mid\alpha_r(s)$ be the
  invariant factors of $P(s)$ and 
let $A$, $\hat \ba$ and  $\hat \bb$ be as in (\ref{eqbigAioicmirmi}), (\ref{eqaioicmirmi}) and (\ref{eqbioicmirmi}), respectively.
From (\ref{sumvu}) we obtain $A=0$, hence (\ref{eqAleq}) holds and from (\ref{eqcdioicmi}) we derive (\ref{eqvusAioi}). 
 Moreover, again from (\ref{sumvu}) we get $\hat \ba =\hat \ba'$ and $\hat \bb =\hat \bb'$,
hence (\ref{eqcmimajhatap}) and  (\ref{eqcmimajpolhatsAioi}) are equivalent 
 and   $\bv=\bu \cup \hat \bb$. By Lemma \ref{lemmacup}, (\ref{eqrmimajpolhatsAioi}) holds.
By Theorem \ref{theoprescrioirmicmipol} and Remark \ref{remacdecr}-\ref{remitac}, the sufficiency of  conditions
(\ref{eqinterpolioi}),
(\ref{eqcdioicmi}) and
 (\ref{eqcmimajhatap}) 
 follows.
\hfill $\Box$

\medskip

\noindent
{\bf Proof of Theorem \ref{corprescrioirmipol}.} 
Assume that there exists a polynomial matrix  $W(s)\in \efe[s]^{z\times n}$  such that
 $\rank \left(\begin{bmatrix}P(s)\\W(s)\end{bmatrix}\right) =r+x$
and $\begin{bmatrix}P(s)\\W(s)\end{bmatrix}$  has 
 $q_1\leq\dots \leq q_{r+x} $ as  invariant orders at $\infty$, 
 and
$\bv=(v_1,\dots, v_{m+z-r-x})$ as row minimal indices. Let
 $\bd=(d_1,\dots, d_{n-r-x})$ be its column minimal indices. By Theorem \ref{theoprescrioirmicmipol},
   (\ref{eqinterpolioi}) and (\ref{eqAleq})--(\ref{eqrmimajpolhatsAioi}) hold where
$A$, $\hat \ba$ and $\hat \bb$ are defined in (\ref{eqbigAioicmirmi}), (\ref{eqaioicmirmi}) and (\ref{eqbioicmirmi}), respectively.

   From (\ref{eqcmimajpolhatsAioi}) we have $d_i\geq c_{i+x}$, $1\leq i \leq n-r-x$; hence $A'=A+\sum_{i=1}^{n-r-x}c_{i+x}-\sum_{i=1}^{n-r-x}d_i\leq A$.
From (\ref{eqAleq}) and (\ref{eqvusAioi})  we obtain (\ref{eqApleq}) and (\ref{eqvusApioi}), respectively.  Since $A'\leq A$, we get 
  $$\sum_{i=1}^j\hat b_i\leq 
  \sum_{i=1}^j\hat b'_i,\quad 1\leq j \leq z-x, 
$$ and from (\ref{eqinterpolioi}),    (\ref{eqAleq}) and (\ref{eqApleq}) we obtain
$$
  \sum_{i=1}^{z-x}\hat b_i=\sum_{i=1}^{m+z-r-x}v_i-\sum_{i=1}^{m-r}u_i=\sum_{i=1}^{z-x}\hat b'_i,
$$
hence $\hat \bb\prec \hat \bb'$,  and from (\ref{eqrmimajpolhatsAioi}) 
and Remark \ref{gmajcp}  we obtain (\ref{eqrmimajpolhatsAioicmi}).

 From (\ref{eqinterpolioi}), we have
$$
\sum_{i=1}^x \hat a'_i =\sum_{i=1}^{m+z-r-x}v_i-\sum_{i=1}^{m-x}u_i+\sum_{i=1}^{r+x}q_i-\sum_{i=1}^{r}p_i-A'=\sum_{i=1}^{x}c_i.
$$
Let $j\in \{1, \dots, x\}$.
Thus,
$$
  \sum_{i=1}^j \hat a_i=
  \sum_{i=1}^j \hat a'_i+A'-A=
\sum_{i=1}^j \hat a'_i+\sum_{i=1}^{n-r-x}c_{i+x}-\sum_{i=1}^{n-r-x}d_i.
$$
Let 
$\hat h_j=\min\{i\; : \; d_{i-j+1}<c_i\}$. Then   $j \leq \hat h_j \leq n-r-x+j$. From (\ref{eqcmimajpolhatsAioi}),
$$
\begin{array}{rl}
  \sum_{i=1}^{\hat h_j}c_i\leq &\sum_{i=1}^{\hat h_j-j}d_i+ \sum_{i=1}^j \hat a_i=\sum_{i=1}^j \hat a'_i+\sum_{i=1}^{n-r-x}c_{i+x}-\sum_{i=\hat h_j-j+1}^{n-r-x}d_i\\
  = & \sum_{i=1}^j \hat a'_i+\sum_{i=1}^{\hat h_j-j}c_{i+x}+\sum_{i=\hat h_j-j+1}^{n-r-x}(c_{i+x}-d_i)
  \leq\sum_{i=1}^j \hat a'_i+\sum_{i=1}^{\hat h_j-j}c_{i+x};
\end{array}
$$
hence
$$\begin{array}{rl}
  \sum_{i=1}^{j}c_i\leq &\sum_{i=1}^j \hat a'_i+\sum_{i=1}^{\hat h_j-j}c_{i+x}-\sum_{i=j+1}^{\hat h_j}c_{i}\\=& \sum_{i=1}^j \hat a'_i+\sum_{i=1}^{\hat h_j-j}(c_{i+x}-c_{i+j})\leq \sum_{i=1}^j \hat a'_i.
\end{array}
$$
Therefore, (\ref{eqcmimajpolhatsAioirmi}) holds.

Conversely, let us assume that
(\ref{eqinterpolioi}) and (\ref{eqApleq})-(\ref{eqrmimajpolhatsAioicmi})
hold. 
Define $d_i=c_{i+x}$, $1\leq i\leq n-r-x$, and $\bd=(d_1, \dots, d_{n-r-x})$.
Let $A$, $\hat \ba$ and $\hat \bb$ be as in 
(\ref{eqbigAioicmirmi}), (\ref{eqaioicmirmi}) and (\ref{eqbioicmirmi}), respectively.
Then $A= A'$, $\hat \ba=\hat \ba'$  and $\hat \bb=\hat\bb'$; hence  (\ref{eqApleq}),
(\ref{eqvusApioi}) and (\ref{eqrmimajpolhatsAioicmi})
 are equivalent to  (\ref{eqAleq}),
(\ref{eqvusAioi}) and (\ref{eqrmimajpolhatsAioi}), respectively.
We have $d_i\geq c_{i+x}$ for $1\leq i \leq n-r-x$, and 
from (\ref{eqinterpolioi}) and (\ref{eqbigAioicmirmi}), 
$\sum_{i=1}^{x}\hat a_i=
\sum_{i=1}^{n-r}c_i - \sum_{i=1}^{n-r-x}d_i$.
Let $j\in \{1, \dots, x\}$ and 
$\hat h_j=\min\{i\; : \; d_{i-j+1}<c_i\}$. 
Then $j \leq \hat h_j \leq n-r-x+j$. 
If
 $j\leq i \leq \hat h_j-1$,  then 
$c_{i+x-j+1}=d_{i-j+1}\geq c_i$; hence
$
c_{i+x-j+1}=c_i$, $j\leq i \leq \hat h_j-1.
$
Thus, 
$$
\begin{array}{rl}
\sum_{i=1}^{\hat h_j}c_i=&\sum_{i=1}^{j-1}c_i+\sum_{i=j}^{\hat h_j-1}c_{i+x-j+1}+c_{\hat h_j}=
\sum_{i=1}^{j-1}c_i+c_{\hat h_j}+\sum_{i=1}^{\hat h_j-j}c_{i+x}\\=&
\sum_{i=1}^{j-1}c_i+c_{\hat h_j}+\sum_{i=1}^{\hat h_j-j}d_{i}\leq \sum_{i=1}^{j}c_i+\sum_{i=1}^{\hat h_j-j}d_{i}.
\end{array}
$$
From (\ref{eqcmimajpolhatsAioirmi}), $\sum_{i=1}^{j}c_i\leq \sum_{i=1}^{j}\hat a'_i$; hence 
$
\sum_{i=1}^{\hat h_j}c_i
\leq \sum_{i=1}^{\hat h_j-j}d_{i}+\sum_{i=1}^{j}\hat a_i.
$
Therefore, (\ref{eqcmimajpolhatsAioi}) holds.
By Theorem \ref{theoprescrioirmicmipol}, the sufficiency of conditions
(\ref{eqinterpolioi}) and (\ref{eqApleq})-(\ref{eqrmimajpolhatsAioicmi})
  follows.
\hfill $\Box$  

\bibliographystyle{acm}
\bibliography{references}

\end{document}